\documentclass[11pt,reqno]{amsart}
\usepackage{graphicx, amssymb, color,pdfsync}
\usepackage[all,cmtip]{xy}
\numberwithin{equation}{section}
\usepackage{mathrsfs}
\usepackage{hyperref}
\newenvironment{psmallmatrix}
  {\left(\begin{smallmatrix}}
  {\end{smallmatrix}\right)}

\hypersetup{
    colorlinks=true,       
    linkcolor=blue,          
    citecolor=blue,        
    filecolor=blue,      
    urlcolor=blue           
}
\usepackage{tikz}
\usetikzlibrary{matrix,arrows}

\usepackage[switch]{lineno}
\usepackage{yfonts}
\advance\hoffset by -0.9 truecm

\begin{document}
\newcommand{\s}{\vspace{0.2cm}}

\newtheorem{theo}{Theorem}
\newtheorem{prop}{Proposition}
\newtheorem{coro}{Corollary}
\newtheorem{lemm}{Lemma}
\newtheorem{claim}{Claim}
\newtheorem{example}{Example}
\theoremstyle{remark}
\newtheorem*{rema}{\it Remark}
\newtheorem*{defi}{\it Definition}
\newtheorem*{theo*}{\bf Theorem}
\newtheorem*{coro*}{Corollary}

\title[Riemann surfaces of genus $2(p-1)$ with $4p$ automorphisms]{The locus of Riemann surfaces of genus $2(p-1)$ \\with $4p$ automorphisms}
\date{}

\author{Angel Carocca}
\address{Departamento de Matem\'atica y Estad\'istica, Universidad de La Frontera, Avenida Francisco Salazar 01145, Temuco, Chile.}
\email{angel.carocca@ufrontera.cl}

\author{Sebasti\'an Reyes-Carocca}
\address{Departamento de Matem\'aticas, Facultad de Ciencias, Universidad de Chile, Las Palmeras 3425, Santiago, Chile}
\email{sebastianreyes.c@uchile.cl}

\thanks{Partially supported by Fondecyt Grants 1240181, 1220099 and 1230708}

\begin{abstract} Let $p \geqslant 5$ be a prime number. In this article we provide a complete and explicit description of the locus formed by the compact Riemann surfaces of genus $2(p-1)$  that are endowed with a group of automorphisms of order $4p$. In addition, we provide isogeny decompositions of the corresponding Jacobian varieties and study if the most symmetric ones admit complex and real multiplication. 
\end{abstract}
\maketitle
\thispagestyle{empty}

\section{Introduction}

The investigation of automorphisms of compact Riemann surfaces (or, equivalently, complex projective algebraic curves and their function fields) and their Jacobian varieties is a classical research field in complex and algebraic geometry. The foundations of the area go back to the nineteenth century, with seminal results given by Riemann, Klein and Jacobi, among others. 

\s

A key result in the area is the uniformisation theorem (see, for instance, \cite[Chapter IV]{FK}). This important theorem states that if a compact Riemann surface is not isomorphic to the projective line nor an elliptic curve, then its universal cover is isomorphic to the upper-half plane $\mathbb{H}$. As a consequence, all such Riemann surfaces can be represented as an orbit space of the form $\mathbb{H}/\Gamma,$ where $\Gamma$ is a Fuchsian group (a discrete subgroup of isometries of $\mathbb{H}$).

\s

Let $X$ be a compact Riemann surface of genus $g \geqslant 2$. A well-known result due to Schwarz says that the full automorphism group of $X$ is finite \cite{Sch1879}. Soon later, Hurwitz by applying his famous ramification formula \cite{Hu} proved that $$|\mbox{Aut}(X)|\leqslant 84(g-1),$$making sense of the problem of determining all groups that act in some given genus. On the other hand, Greenberg in \cite{Gre} proved that each  finite group can be realised as the automorphism group of some compact Riemann surface, and therefore of some Jacobian variety.

\s

Given a complex projective algebraic curve of genus $g \geqslant 2$ through polynomial equations, it is rather difficult to get information about its automorphisms. To deal with this problem, the uniformisation theorem and the theory of Fuchsian groups provide a powerful toolkit: they allow us to describe and construct group action on Riemann surfaces, as such groups arise as quotients of Fuchsian groups. 

\s

Articles aimed at studying group actions on Riemann surfaces and Jacobian varieties can be found in the literature in plentiful supply. We refer to \cite{B22}, \cite{BuCi} and \cite{Rubi14} (and to the references therein) for an up-to-date treatment of this topic.

\s

Let $\mathscr{M}_g$ denote the moduli space of compact Riemann surfaces of genus $g.$ It is known that  $\mathscr{M}_g$ has a structure of a complex analytic space of dimension $3g-3$, whose singularities agree with the points representing Riemann surfaces with non-trivial automorphisms, that is, $$\mbox{Sing}(\mathscr{M}_g)=\{[X] \in \mathscr{M}_g : \mbox{Aut}(X) \neq \{\mbox{id}\}\} \mbox{ for } g \geqslant 4.$$

The singular locus of $\mathscr{M}_g$ becomes increasingly complicated as the genus grows. For instance, the diversity of automorphism groups as well as the different ways in which such groups can act (which is closely related with the topology of $\mathscr{M}_g$) increase significantly. 

\s

Let $\mathcal{A}_g$ denote the moduli space of principally polarised abelian varieties of dimension $g.$ As in the case of Riemann surfaces, $\mathcal{A}_g$ is endowed with a structure of complex analytic space of dimension $g(g+1)/2$ and $$\mbox{Sing}(\mathcal{A}_g)=\{[A] \in \mathcal{A}_g : \mbox{Aut}(A) \neq \{\pm \mbox{id}\}\} \mbox{ for } g \geqslant 3.$$

Let $X \in \mathscr{M}_g$ and let $JX$ be its  Jacobian variety. It is well-known that $JX \in \mathcal{A}_g$ is irreducible and that the so-called Torelli map $$\mathscr{M}_g \to \mathcal{A}_g \, \mbox{ defined by } \,  [X] \mapsto [JX]$$is injective; see, for instance, \cite[Chapter III]{FK}. In other words, compact Riemann surfaces are determined by their Jacobian varieties.


\s

To study the largely unknown topology of the singular locus of $\mathscr{M}_g$, the loci formed by Riemann surfaces endowed with a group of automorphisms of a prescribed order $$\mathscr{M}_g^{n}=\{[X] \in \mathscr{M}_g : X \mbox{ has a group of automorphisms of order }n\}$$ have been widely studied. Among such prescribed orders, the ones that are linear in the genus,  $$ag+b \mbox{ where $a, b$ are rationals,}$$have played a significant role. A remarkable case is $4g+2$, which was studied by Wiman in \cite{Wi}. He proved that this number is the largest order among the cyclic groups of automorphisms of Riemann surfaces of genus $g$ (see also \cite{Harvey1}). Wiman also showed that such a bound is sharp, by exhibiting the Riemann surface $W_g$ represented by the curve \begin{equation*} \label{eWiman}y^2=x^{2g+1}-1 \mbox{ and the automorphism } (x,y) \mapsto (\omega_{2g+1}x, -y).\end{equation*}Kulkarni went further by proving in \cite{K1} that $W_g$ is uniquely determined by this property.

\s

The case $8g+8$ has also been completely described, yielding the so-called Accola-Maclachlan Riemann surfaces; see \cite{A68}, \cite{Mac} and also the survey \cite{Sing14}. In a very similar line of research, the case $4g$ was considered in \cite{BCI17}, and the case $\lambda(g-1)$, where $\lambda \geqslant 3$ is an integer and $g-1$ is prime, was considered first in \cite{BJ05} and completed recently in \cite{IJRC21}. We also refer to  \cite{CRC},  \cite{CY92} and \cite{CI18} for more cases. The analogous problem for Klein surfaces has been also considered. See, for instance, \cite{BE88}, \cite{M77} and the recent article \cite{AC23} (and the references therein).

\s

Let $p \geqslant 5$ be a prime number. In this article, we provide a complete and explicit description of the sublocus of $\mathscr{M}_{2(p-1)}$ formed by those compact Riemann surfaces that are endowed with a group of automorphisms of order $4p$. Throughout the article, this locus will be denoted by  $$\mathscr{M}_{2(p-1)}^{4p}=\{[X] \in \mathscr{M}_{2(p-1)} : X \mbox{ has a group of automorphisms of order }4p\}.$$ 

In our case, the number of automorphisms that the Riemann surfaces admit is $2g+4$. It is worth remarking that, in contrast to most of the cases already mentioned, $2g+4$ is relatively small compared to the genus.
\s

The results of this paper --that will be stated in Section \S \ref{statement}-- can be succinctly summarised as follows.

\s

{\bf 1.} We study the topology of $\mathscr{M}_{2(p-1)}^{4p}$ by taking advantage of the equisymmetric stratification of $\mathscr{M}_g$ introduced in \cite{B90}. We prove that this locus has complex dimension one, and is formed by two disjoint complex one-dimensional families, together with $$\tfrac{2p-5}{3} \mbox{ if } p \equiv 1 \mbox{ mod } 3 \,\,  \mbox{ and } \,\, \tfrac{2p-7}{3} \mbox{ if }p \equiv 2 \mbox{ mod } 3$$ pairwise nonisomorphic quasiplatonic Riemann surfaces (that do not belong to the families).

\s

{\bf 2.} By the way, we obtain that there is no compact Riemann surface in $\mathscr{M}_{2(p-1)}$ endowed with more than $12p$ automorphisms.

\s

{\bf 3.} For each member $X$ of $\mathscr{M}_{2(p-1)}^{4p}$ we determine its full automorphism group and describe the way it acts. Also, we prove that $X$ is cyclic $p$-gonal and non-hyperelliptic, and we describe it as an explicit affine plane algebraic curve.

\s

{\bf 4.} We study the corresponding Jacobian varieties by providing isogeny decompositions of them. Further, we provide the Poincar\'e isogeny decomposition of a distinguished member of the locus, and consider the problem of studying whether the most symmetric members admit complex and real multiplication. \s

This article is organised as follows. In Section \S\ref{preli} we introduce some notation and briefly review the basic preliminaries. Section \S\ref{statement} is devoted to the statement of the results. The proofs will be given in Section \S\ref{proof1}.
 
\section{Preliminaries} \label{preli}

Let $X$ be a Riemann surface of genus $g \geqslant 2$ and let $G$ be a finite group acting on $X$ $$\psi: G \to \mbox{Aut}(X).$$ Consider the quotient map $\pi_G : X \to X/G=X_G$ given by the action of $G \cong \psi(G)$ on $X$, and assume that it ramifies over $r$ values $y_1, \ldots, y_r$. Let $h$ denote the genus of $X/G.$ The {\it signature} of the action of $G$ is the tuple $$s=(h; m_1, \ldots, m_r) \in \mathbb{Z}^+_0 \times (\mathbb{Z}_{\geqslant 2})^r$$where  $m_i$ is the order of the $G$-stabiliser of $x_i \in X$ where $\pi_G(x_i)=y_i$.

\s

Let $\Gamma$ be a Fuchsian group such that $X \cong \mathbb{H}/\Gamma.$ Riemann existence theorem ensures that $G$ acts on $X$ if and only if there is a Fuchsian group $\Delta$ and there is an epimorphism  $$\theta : \Delta \to G \mbox{ such that } \mbox{ker}(\theta)=\Gamma.$$ The genus $g$ of $X$ is related to $s$ by the {\it Riemann-Hurwitz formula} $$2(g-1)=|G|(2(h-1)+\sum_{j=1}^r(1-1/m_j)).$$

The epimorphism $\theta$ represents the action and is called a {\it surface-kernel epimorphism} (ske, for short). We usually identify $\theta$ with the tuple or {\it generating vector} $$\theta=(\theta(\alpha_1), \ldots, \theta(\alpha_h), \theta(\beta_1), \ldots, \theta(\beta_h), \theta(\gamma_1), \ldots, \theta(\gamma_r)) \in G^{2h+l}.$$

As $X/G \cong \mathbb{H}/\Delta$, we say that the signature of $\Delta$ is $s$, and write $\Delta=\Delta_s.$ It is classically known that $\Delta_s$ admits a canonical presentation is terms of $2h$ hyperbolic generators $\alpha_1, \beta_1, \ldots, \alpha_h, \beta_h$ and $r$ elliptic generators $\gamma_1, \ldots, \gamma_r$ with relations $$\prod_{i=1}^h\alpha_i\beta_i\alpha_i^{-1}\beta_i^{-1}\prod_{j=1}^r\gamma_j=\gamma_1^{m_1}=\cdots=\gamma_r^{m_r}=1.$$The elements of finite order of $\Delta_s$ are those that are conjugate to powers of $\gamma_1, \ldots, \gamma_r.$

\s

Let $G'$ be a supergroup of $G.$ The action of $G$ on $X$ is said to {\it extend} to an action of $G'$ if the following three statements hold.\begin{enumerate}
\item There is a Fuchsian group $\Delta'$ with $\Delta \leqslant \Delta'$.
\item There is a surface-kernel  epimorphism $$\Theta: \Delta' \to G' \, \, \mbox{ with }  \, \, \Theta|_{\Delta}=\theta.$$
\item The Teichm\"{u}ller spaces of $\Delta$ and $\Delta'$ have the same dimension.
\end{enumerate} 
An action is called {\it maximal} if it cannot be extended. We refer to  \cite{Sing72} for the list of signatures that might admit a non-maximal action. See also \cite{Bu03} and \cite{Ries93}.  

\s

Let $Y$ be a  Riemann surface of genus $g \geqslant 2$ endowed with an action $\psi':G \to \mbox{Aut}(Y)$. If there exists an orientation-preserving homeomorphism  \begin{equation} \label{lapiz22}\phi: X \to Y \mbox{ such that } \phi \psi(G) \phi^{-1}=\psi'(G)\end{equation}then it is said that the actions of $G$ are {\it topologically equivalent}. If we write $$X/{G} \cong \mathbb{H}/{\Delta} \, \mbox{ and } \, Y/G \cong \mathbb{H}/{\Delta}'$$then each orientation-preserving  homeomorphism $\phi$ as in \eqref{lapiz22}  induces  a group isomorphism $\phi^*: \Delta \to \Delta'$. Hence, we may assume $\Delta=\Delta'$. We denote the subgroup of $\mbox{Aut}(\Delta)$ consisting of the automorphisms $\phi^*$ by $\mathscr{B}$. Following \cite{B90}, the skes $$\theta : \Delta \to G \, \mbox{ and } \, \theta': \Delta \to G$$representing two actions of $G$  are topologically equivalent if and only if $$\mbox{there exist } a \in \mbox{Aut}(G) \mbox{ and } \phi^* \in \mathscr{B} \mbox{ such that } \theta' = a\circ\theta \circ \phi^*.$$ In such a case, we write $\theta \sim \theta'.$ We refer to Nielsen \cite{Nielsen}, Harvey \cite{H71} and Gilman \cite{Gilman} as sources for the characterisation of topological actions by certain purely algebraic data. 
\s

A {\it family of Riemann surfaces} of genus $g$ is the locus $$\mathcal{F}=\mathcal{F}(G, s) \subset \mbox{Sing}(\mathscr{M}_g)$$ formed by the points representing isomorphism classes of Riemann surfaces endowed with an action of $G$ with signature $s$. According to \cite{B90}, the interior of the family $\mathcal{F}$, if non-empty, consists of those Riemann surfaces whose (full) automorphism group is isomorphic to $G$ and is formed by finitely many components, each one of them being a closed irreducible (non necessarily smooth) subvariety of  $\mathscr{M}_g$ (see also \cite{GG92}). 

\s

We refer to the articles \cite{BCI14} and \cite{CI10}  for the case of low genus, and to \cite{BRT23}, \cite{Pign} and  \cite{Jen17} for algorithms to find concrete examples of subvarieties of this kind. 

\s

{\it Notation.} We denote the cyclic group of order $n$ by $\mathbb{Z}_n$, the dihedral group of order $2n$ by $\mathbf{D}_n$ and the alternating group of order $12$ by $\mathbf{A}_4.$ We write $\omega_p=\mbox{exp}(2\pi i/p).$
\section{Statement of the Results}\label{statement}

 Let $p \geqslant 5$ be a prime number.  Throughout this paper, we write $g_p:=2(p-1)$ and denote by $X_p$ and $Y_p$ the  Riemann surfaces of genus $g_p$ represented by the curves $$y^p=x^{4}(x^4-1)^{p-2} \,\,\, \mbox{ and } \,\,\,  y^p=(x^3-1)(x^3+1)^{p-1}$$respectively. Also, for each $p \equiv 1 \mbox{ mod }3$, we denote by $Z_p$  the  Riemann surface of genus $g_p$ represented by the  curve  $$y^p=x^2(x^2-1)^{2r}(x^2+1)^{2r^2}$$where $r$ is a primitive third root of unity in the field of $p$ elements.
 
 \s

%
%
%
%
%
%
%

The following result shows that the aforementioned Riemann surfaces enjoy the property of being determined by the number of automorphisms they admit.

\begin{theo}\label{teo8p12p} Let $p \geqslant 5$ be a prime number. 

\s

{\bf 1.} If $X$ is a compact Riemann surface of genus $g_p$ endowed with a group of automorphisms of order $8p$ then $X \cong X_p.$ 

\s

{\bf 2.} If $X$ is a compact Riemann surface of genus $g_p$ endowed with a group of automorphisms of order $12p$, then
\begin{enumerate}
\item[(a)] if $p \equiv 1 \mbox{ mod } 3$ then $X \cong Y_p$ or $X \cong Z_p,$ and
\item[(b)] if $p \equiv 2 \mbox{ mod } 3$ then $X \cong Y_p$.
\end{enumerate}
In addition, 
$$\mbox{Aut}(X_p) \cong \mathbb{Z}_p \times \mathbf{D}_4, \,\, \mbox{Aut}(Y_p) \cong \mathbf{D}_3 \times \mathbf{D}_p, \,\,\mbox{Aut}(Z_p)\cong \mathbb{Z}_p \rtimes_3 \mathbf{A}_4$$and these groups act on the corresponding Riemann surfaces with signatures $$(0; 2,2p,4p),\,\, (0; 2,6,2p) \, \mbox{ and } \, (0;3,3,2p).$$ 
\end{theo}

Observe that $X_p, Y_p, Z_p$ belong to $\mathscr{M}_{2(p-1)}^{4p}$ as they have a group of automorphisms of order $4p.$ Indeed, if $X \in \mathscr{M}_{2(p-1)}$ then $$X \in \mathscr{M}_{2(p-1)}^{4p} \mbox{ if and only if } |\mbox{Aut}(X)|=4 \lambda p \mbox{ for some }\lambda \geqslant 1.$$ The following result shows that there are only three possibilities for $\lambda$.

\begin{theo}\label{theobig} Let $p \geqslant 5$ be a prime number. There is no compact Riemann surface of genus $g_p$ endowed with a group of automorphisms of order $4\lambda p$ for  $\lambda \geqslant 4.$
\end{theo}

The theorem above coupled with Theorem \ref{teo8p12p} allow us to focus our attention on Riemann surfaces endowed with a group of automorphisms of order $4p.$ The following lemma provides the potential groups and potential signatures in such a case.

\begin{prop}\label{lema1} Let $p \geqslant 5$ be a prime number and let $G$ be a group of order $4p$ acting in genus $g_p$. The signature of the action  is 
$$s_1:=(0; 2,2,p,2p), \, s_2:=(0; 2p, 2p, 2p) \, \mbox{ or }\, s_3:=(0; p, 4p, 4p).$$In addition, the following statements hold.
\begin{enumerate}
\item $G$ acts with signature $s_1$ if and only if $G$ is isomorphic to $\mathbb{Z}_p \times \mathbb{Z}_2^2$ or $\mathbf{D}_{2p}.$ 
\item  $G$ acts  with signature $s_2$ if and only if $G$ is isomorphic to $\mathbb{Z}_p \times \mathbb{Z}_2^2$.
\item  $G$ acts with signature $s_3$ if and only if $G$ is cyclic.
\end{enumerate}
\end{prop}

We anticipate that the potential groups and signatures in the proposition above are realised. The following result describes the Riemann surfaces (and the locus in the moduli space formed by them) corresponding to case (1) in the proposition above.

\begin{theo} \label{teoA}
Let $p \geqslant 5$ be a prime number. 

\s

{\bf 1.} The Riemann surfaces of genus $g_p$ endowed with a group of automorphisms $G$ isomorphic to $\mathbb{Z}_p \times \mathbb{Z}_2^2$ acting with signature $s_1$ form a complex one-dimensional family.

\s

Moreover, if we denote such a family by $\mathcal{F}_1$ then the following statements hold.

\begin{enumerate}

\item $\mathcal{F}_1$ is an irreducible subvariety of $\mathscr{M}_{2(p-1)}.$

\s

\item The members of $\mathcal{F}_1$  are $p$-gonal and algebraically represented by the curves $$y^p=x^{4}[(x^2-1)(x^2-t)]^{p-2} \mbox{ where }t \in \mathbb{C}-\{0,1\}.$$

\item In the model above, the group $G$ is generated by the transformations $$(x,y) \mapsto (x,\omega_p y), \,\, (x,y) \mapsto (-x,y)\, \mbox{ and } \,(x,y) \mapsto(\tfrac{\sqrt{t}}{x}, \tfrac{ty}{x^4}).$$

\item The Riemann surface $X_p$ belongs to $\mathcal{F}_1.$ 

\s

\item If $X \in \mathcal{F}_1$ is not isomorphic to $X_p$ then its full automorphism group  is $G$.

\end{enumerate}

\s

{\bf 2.} The Riemann surfaces of genus $g_p$ endowed with a group of automorphisms isomorphic to $\mathbf{D}_{2p}$ acting with signature $s_1$ form a  complex one-dimensional  family.

\s

Moreover, if we denote such a family by $\mathcal{F}_2$ then the following statements hold.

\begin{enumerate}

\item $\mathcal{F}_2$ consists of at most $\tfrac{p-1}{2}$ irreducible components.

\s

\item The members of $\mathcal{F}_2$ are $p$-gonal and algebraically represented by the curves $$y^p=(x-1)(x-t)^k(x-\tfrac{1}{t})^k(x+1)^{p-1}(x+t)^{p-k}(x+\tfrac{1}{t})^{p-k}$$where $t \in \mathbb{C}-\{0, \pm1, \pm \sqrt{-1}\}$ and $k \in \{1, \ldots, \tfrac{p-1}{2}\}.$

\s

 \item In the model above, the group $G$ is generated by the transformations $$(x,y) \mapsto (\tfrac{1}{x},-\tfrac{\omega_p y}{x^3}) \mbox{ and } (x,y) \mapsto (-x, \tfrac{1}{y}(x^2-1)(x^2-t^2)(x^2-\tfrac{1}{t^2}))$$

\item The Riemann surface $Y_p$ belongs to $\mathcal{F}_2.$ 

\s

\item If $X \in \mathcal{F}_2$ is not isomorphic to $Y_p$ then its full automorphism group is $G$.

\end{enumerate}
\end{theo}

\s

If $\mathcal{F} \subset \mathscr{M}_g$ is a family of  Riemann surfaces endowed with a group of automorphisms $G$ acting with a maximal signature (as we shall see that is the case of $\mathcal{F}_1$ and $\mathcal{F}_2$), then we write $$\partial(\mathcal{F}):=\{X \in \mathcal{F}: G \mbox{ is properly contained in }\mbox{Aut}(X)\}.$$ 
The maximality of the signature implies that $\partial(\mathcal{F})$ has positive codimension in $\mathcal{F}.$ As a consequence of the theorem above coupled with Theorem \ref{theobig}, one obtains the following corollary.
\begin{coro}\label{coro1}
$\mathcal{F}_1 \cap \mathcal{F}_2=\emptyset, \partial(\mathcal{F}_1)=\{X_p\}$ and $\partial(\mathcal{F}_2)=\{Y_p\}.$
\end{coro}

We now proceed to study the Riemann surfaces corresponding to case (2) in Proposition \ref{lema1}. For each $j \in \{1, \ldots, p-2\}$, let $C_j$ be the  Riemann surface of genus $g_p$ represented by $$y^p=x^{2j+2}(x^2-1)^{p-2}(x^2+1)^{p-2j}$$where the exponents are taken modulo $p$ and lie in $\{1, \ldots, p-1\}.$ Note that $C_j$ has a group of automorphisms $G$ isomorphic to $\mathbb{Z}_p \times \mathbb{Z}_2^2$  generated by  $$ (x,y) \mapsto (x, \omega_p y), \, (x,y) \mapsto (-x, y) \mbox{ and } (x,y) \mapsto (\tfrac{1}{x}, -\tfrac{y}{x^4})$$In particular, each $C_j$ belongs to $\mathscr{M}_{2(p-1)}^{4p}.$

\begin{theo} \label{teoD} Let $p \geqslant 5$ be a prime number.  If $X$ is a Riemann surface of genus $g_p$ endowed with a group of automorphisms $G$ isomorphic to $\mathbb{Z}_p \times \mathbb{Z}_2^2$ acting with signature $s_2$ then $$X \cong C_j \mbox{ for some } j \in \{1, \ldots, p-2\}.$$
\s
Moreover, the following statements hold.
\begin{enumerate}
\item  $C_{j_1}$ is isomorphic to $C_{j_2}$ if and only if $j_1=\sigma(j_2)$ where $\sigma \in \langle z \mapsto \tfrac{1}{z}, z \mapsto \tfrac{-1}{1+z}\rangle.$

\s

\item $C_{1}$ is isomorphic to $X_p.$ 

\s

\item If $p \equiv 1 \mbox{ mod }3$ and $2 \leqslant r \leqslant p-2$ satisfies $r^3 \equiv 1 \mbox{ mod }p,$ then $C_r \cong Z_p.$

\s

\item If $C_j$ is not isomorphic to $X_p$ nor $Z_p$ then the full automorphism group of $C_j$  is $G$.

\s

\item There are exactly$$\tfrac{p+5}{6} \mbox{ if } p \equiv 1 \mbox{ mod } 3 \,\,  \mbox{ and } \,\, \tfrac{p+1}{6} \mbox{ if }p \equiv 2 \mbox{ mod } 3$$pairwise nonisomorphic Riemann surfaces among $C_1, \ldots, C_{p-2}.$

\end{enumerate}
\end{theo}

An interesting consequence of the theorems above together with the results of \cite{GG92} is the fact that $\mathcal{F}_1$ and each irreducible component of $\mathcal{F}_2$ do not have non-normal points. Hence, the following result is obtained.

\begin{coro}\label{normal}
The family $\mathcal{F}_1$ and each irreducible component of $\mathcal{F}_2$ are subvarieties of $\mathscr{M}_{2(p-1)}$ that are smooth (and hence, they themselves are Riemann surfaces).
\end{coro}

To complete the classification, we now consider the Riemann surfaces corresponding to case (3) in Proposition \ref{lema1}. For each $j \in \{1, \ldots, \frac{p-1}{2}\}$, let $S_j$ be the  Riemann surface of genus $g_p$ represented by the algebraic curve $$y^p=x^4(x^4-1)^s \mbox{ where } s \in \{1, \ldots, p-1\} \mbox{ satisfies } sj \equiv 1 \mbox{ mod }p.$$As $(x,y) \mapsto (ix,\omega_px)$
is an automorphism of $S_j$ of order $4p$, each $S_j$ belongs to $\mathscr{M}_{2(p-1)}^{4p}.$

\begin{theo} \label{theoC} Let $p \geqslant 5$ be a prime number. If $X$ is a Riemann surface of genus $g_p$ endowed with a cyclic group $G$ of automorphisms of order $4p$ then $$X \cong S_j \mbox{ for some }j \in \{1, \ldots, \tfrac{p-1}{2}\}.$$

Moreover, the following statements hold.

\begin{enumerate}

\item  $S_{j_1}$ is isomorphic to $S_{j_2}$ if and only if $j_1=j_2.$

 \s
 
\item $S_{\frac{p-1}{2}}$ is isomorphic to $X_p.$

\s

\item If $j \neq \frac{p-1}{2}$ then the full automorphism group of $S_j$  is $G$.
\end{enumerate}
\end{theo}

\s

The theorems above can be summarised in the following result, as anticipated in the introduction.
\begin{theo}
Let $p \geqslant 5$ be a prime number. Then  $\mathscr{M}_{2(p-1)}^{4p}$ has complex dimension one, and is formed by two disjoint complex one-dimensional families together with $$\tfrac{2p-5}{3} \mbox{ if } p \equiv 1 \mbox{ mod } 3 \,\,  \mbox{ and } \,\, \tfrac{2p-7}{3} \mbox{ if }p \equiv 2 \mbox{ mod } 3$$ pairwise nonisomorphic quasiplatonic Riemann surfaces (that do not belong to the families).

\end{theo}

According to the main result of \cite{Wootton05}, if a compact Riemann surface is simultaneously hyperelliptic and $p$-gonal then its full automorphism group is isomorphic to $\mathbb{Z}_{2p}, \mathbb{Z}_p \rtimes \mathbf{D}_4$ or $\mathbf{D}_{2p}$, with this latter group acting with signature $(0; 2,2,2,p)$.

As our Riemann surfaces are $p$-gonal, we obtain the following result.
\begin{coro}
$\mathscr{M}_{2(p-1)}^{4p}$ does not intersect the hyperelliptic locus of $\mathscr{M}_{2(p-1)}.$
\end{coro}

We now proceed to provide isogeny decomposition of  $JX$ for each  
 $X  \in \mathscr{M}_{2(p-1)}^{4p}$.
 
 \begin{theo}\label{jacofam}
 Let $p \geqslant 5$ be prime. For each $X \in \mathscr{M}_{2(p-1)}^{4p}$ consider the group of automorphisms $G$ as in Theorems \ref{teoA}, \ref{teoD} and \ref{theoC}. If $H \leqslant G$ then we write $X_H=X/H.$
 
 \s
 
 {\bf 1.} If $X \in \mathcal{F}_1$ and $\sigma, \tau$ are the two involutions in $G$ with $2p$ fixed points, then $$JX \sim JX_{\langle \sigma \rangle} \times JX_{\langle \tau \rangle} \times JX_{\langle \sigma \tau \rangle}$$where $JX_{\langle \sigma \rangle}$ and $JX_{\langle \tau \rangle}$ have dimension $\tfrac{p-1}{2}$ and $JX_{\langle \sigma \tau \rangle}$ has dimension $p-1.$ 
 
 \s
 
 {\bf 2.} If $X \in \mathcal{F}_2$ and $\iota$ is any non-central involution in $G$ then $$JX \sim JX_{\langle \iota \rangle}^2.$$Further,  $JX_{\langle \iota \rangle}$ is isogenous to the product of two factors of dimension $\tfrac{p-1}{2}.$
 
 \s

 {\bf 3.} If $\tau_1$ and $\tau_2$ are two involutions in $G$ and $H=\langle \tau_1, \tau_2\rangle$ then $$JC_j\sim J(C_j)_{H} \times \prod_{t=1}^3 \mbox{Prym}((C_j)_{H_i} \to (C_j)_{H}) \mbox{ for each }j,$$where $H_1, H_2$ and $H_3$ are the three subgroups of $H$ of order two. Each one of the factors has dimension $\tfrac{p-1}{2}.$
 
 \s
 
 {\bf 4.} If $\iota$ denotes the unique involution in $G$ then $$JS_j \sim J(S_j)_{\langle \iota \rangle} \times P \mbox{ for each }j,$$where each factor has dimension $p-1.$ Further, the former factor is isogenous to the product of two factors of dimension $\tfrac{p-1}{2}.$
 \end{theo}
 
The theorem above implies that $JX_p \cong JC_1$ admits an isogeny decomposition in terms of four non necessarily simple factors that are pairwise non-$K$-isogenous, with $K$ being each one of the three non-conjugate subgroups of $\mbox{Aut}(X_p)$ of order $4p$.  However, the fact that $X_p$ admits extra automorphisms allows us to go further by providing its Poincar\'e isogeny decomposition (that is, a decomposition in terms of simple factors).  
\begin{theo} \label{jaco} Let $p \geqslant 5$ be prime. 
The Poincar\'e isogeny decomposition of  $JX_p$ is $$JX_p \sim J\hat{X}_p^4$$where $\hat{X}_p$ is the hyperelliptic Riemann surface of genus $\tfrac{p-1}{2}$  represented by  the curve$$y^p=x(x-1)^{p-2}.$$
\end{theo}

Let $A$ be abelian variety of dimension $g$ and let $\mbox{End}_{\mathbb{Q}}(A)$ denote its rational endomorphism algebra. It is said that 
\begin{enumerate}
\item $A$ has complex multiplication if  there exists a CM field $K$ such that $$K \subset \mbox{End}_{\mathbb{Q}}(A) \mbox{ and } [K: \mathbb{Q}]=2g.$$
\item $A$ has real multiplication if  there exists a totally real number field $K$ such that $$K \subset \mbox{End}_{\mathbb{Q}}(A) \mbox{ and } [K: \mathbb{Q}]=g.$$
\end{enumerate}

Little is currently understood about the existence of curves and families of curves whose Jacobians are acted upon by large rings of endomorphisms. Consequently, constructing curves with Jacobians exhibiting real or complex multiplication is an area of significant interest. Jacobians with complex multiplication have attracted considerable interest in algebraic geometers and number theorists, as they represent special points in $\mathcal{A}_g$ and enjoy interesting arithmetic properties. For instance, Jacobians with complex multiplication can be defined over number fields \cite{ST61}, providing an interesting relationship with triangle Fuchsian groups and the theory of {\it dessin d'enfants} introduced by Grothendieck (see the survey \cite{JW}). For explicit examples of Jacobians with real and complex multiplication we refer to \cite{CLR11}, \cite{e}, \cite{M91} and \cite{T91}.

\begin{theo}\label{CM}
Let $p \geqslant 5$ be a prime number. 
\begin{enumerate}
\item $JX_p$ has complex multiplication by $\mathbb{Q}(\omega_p)$.

\s

\item $JY_p$ has real multiplication by 
$\mathbb{Q}(\omega_p+\bar{\omega}_p)$.

\s

\item $JZ_p$ is isogenous to the product of two factors with complex multiplication by $\mathbb{Q}(\omega_p)$ and $\mathbb{Q}(\omega_p+\omega_p^r+\omega_p^{r^2})$ respectively, where $r$ is a primitive third root of unity in the field of $p$ elements.
\end{enumerate}
\end{theo}

\section{Proofs}\label{proof1}
\subsection*{Proof of Theorem \ref{teo8p12p}} \mbox{}

{\bf 1.} Assume that $X \in \mathscr{M}_{2(p-1)}$ has a group of automorphisms $G$ of order $8p.$ By the Riemann-Hurwitz formula,  the signature of the action of $G$ on $X$  has the form \begin{equation}\label{rh}(0; m_1, m_2, m_3) \, \mbox{ where } \, \tfrac{1}{m_1}+\tfrac{1}{m_2}+\tfrac{1}{m_3}=\tfrac{1}{2}+\tfrac{3}{4p}.\end{equation}

\s

{\it Claim.} The signature of the action of $G$ on $X$ is $(0; 2,2p,4p).$

\s

Let $v$ be the number of periods $m_i$ that are equal to two. Clearly $v \neq 2,3.$ Assume $v=0$. If each $m_i$ were greater or equal to 6 then equation \eqref{rh} would not be satisfied. Thus, without loss of generality, we have that $m_1=4$ or $m_1=5.$ 

\s

Assume $m_1=4.$ In such a case  equation \eqref{rh} turns into \begin{equation}\label{rh2}\tfrac{1}{m_2}+\tfrac{1}{m_3}=\tfrac{1}{4}+\tfrac{3}{4p} \mbox{ where }m_2, m_3 \geqslant 4.\end{equation}It is easy to see that $m_2 \neq 4,6,7.$ If $m_2=5$ then $m_3=p=5.$ However, there is no Riemann surface of genus $8$ with a group of automorphisms of order 40 acting with signature $(0; 4,5,5).$ See, for instance, \cite{Conder} and \cite{K}. Thus, $m_2, m_3 \geqslant 8$ but this contradicts \eqref{rh2}. 

\s

Assume $m_1=5.$ In such a case $p=5$ and  equation \eqref{rh} turns into \begin{equation}\label{rh3}\tfrac{1}{m_2}+\tfrac{1}{m_3}=\tfrac{9}{20} \mbox{ where }m_2, m_3 \geqslant 4.\end{equation}Note that if $m_2=4$ then $m_3 \notin \mathbb{Z}$. Thus $m_2, m_3 \geqslant 5,$ contradicting \eqref{rh3}.

\s

All the above says that the signature of the action of $G$ on $X$ is $$(0; 2, m_2, m_3) \,\mbox{ where } \, \tfrac{1}{m_2}+\tfrac{1}{m_3}=\tfrac{3}{4p}.$$

We assume that $m_2, m_3 \neq 2p.$ If $m_2 < 2p$ then $m_3>4p$ and hence $m_3=8p.$ In such a case, if the action is represented by a ske $(u,v,w)$ then, as $G$ is cyclic, one has that the order of $v=(wu)^{-1}$ is $8p,$ contradicting the assumption that $m_2 < 2p$. Similarly, if $m_2 > 2p$ then $m_3$ equals to $4,8$ or $p,$ but none of these possibilities can be realised. Hence, $m_2=2p$ and $m_3=4p,$ and the signature of the action is $(0; 2,2p,4p)$, as claimed.

\s

{\it Claim.} $G \cong \mathbb{Z}_p \times \mathbf{D}_4$ and $X \cong X_p.$

\s

Assume that the action of $G$ on $X$ is represented by the ske $$\theta: \Delta_{(0; 2, 2p, 4p)} \to G \mbox{ given by }\theta = (x,y,z)$$and let $N =\langle z \rangle.$ Observe that $x \notin N$, since otherwise $\theta$ is non-surjective. If we write $a:=z^4,  r:=z^p$ and $s:=x$ then one has that $$G \cong \langle a,r,s : a^p=r^4=s^2=1, [a,r]=1, sas=a^u, srs=r^v\rangle \mbox{ where }u, v=\pm 1.$$Note that $G$ cannot be abelian, and thus the case $u=v=1$ must be disregarded.

\s

If $u=-1$ and $v=1$ then $G$ is isomorphic to $\mathbf{D}_p \times \mathbb{Z}_4.$ In such a case, the elements of order $4p$ are $a^ir^{\pm 1}$ whereas the elements of order $2p$ are $a^ir^2$ where $i \in \{1, \ldots, p-1\}.$ Observe that $G$ cannot be generated by one element of order $4p$ and one element of order $2p$, showing that this case is not realised.

\s

If $u=v=-1$ then $G$ is isomorphic to the semidirect product $$\langle a,r \rangle \rtimes \langle s \rangle \cong (\mathbb{Z}_p \times \mathbb{Z}_4) \rtimes \mathbb{Z}_2 \, \mbox{ where } (sg)^2=1 \mbox{ for all } g \in \langle a,r \rangle.$$As the elements of this group of order $2p$ and $4p$ are contained in the normal subgroup $\langle a,r \rangle$, this case is not realised neither. 

\s

We then conclude that necessarily $u=1, v=-1$ and therefore  $G \cong  \mathbb{Z}_p \times \mathbf{D}_4.$

\s

The elements of order two of $G$ are $r^2$ and $sr^k$, the elements of order $2p$ are $a^ir^2$ and $a^isr^k$, whereas the ones of order $4p$ are $a^ir^{\pm 1}$ where $k \in \{1, \ldots, 4\}$ and $i \in \{1, \ldots, p-1\}.$ Then, after considering  an automorphism of $G$, we can assume that $$\theta=(\eta, \eta^{-1}a^{-1}r, ar^{-1}) \mbox{ for some involution }\eta \in G.$$ If $\eta$ were equal to $r^2$ then $\theta$ would not be surjective.  Thus, up to an automorphism of $G$, we can assume $\eta=s$. It follows that $\theta$ is $G$-equivalent to $(s,a^{-1}sr,ar^{-1})$ and hence $X$ is uniquely determined, up to isomorphism. The claim follows from the fact that 
$$(x,y) \mapsto  (x, \omega_py), \,\, (x,y) \mapsto (ix, y) \,\, \mbox{ and } \,\, (x,y)\mapsto (\tfrac{i}{x}, -\tfrac{y}{x^4})$$generate a group of automorphisms of $X_p$ isomorphic to $\mathbb{Z}_p \times \mathbf{D}_4$.

\s

{\it Claim.} $\mbox{Aut}(X) = G.$
\s

Assume that $G$ is strictly contained in the full automorphism group $\hat{G}$ of $X$. Then, as a consequence of the results of \cite{Sing72}, this group has order $24p$, acts on $X$ with signature $(0; 2,3,4p)$ and contains $G$ as a non-normal subgroup. The natural action of $\hat{G}$ on the set of right cosets $\hat{G}/G$ induces a homomorphism  $\Theta :  \hat{G} \to {\bf S}_3.$ Therefore there is a subgroup $L$ of $\mathbf{S}_3$ such that 
$$\hat{G}/K \cong L \leqslant {\bf S}_3 \, \mbox{ where } \, K = \ker(\Theta) = \bigcap_{g \in \hat{G}} G^g  \leqslant G $$
Since $G$ is not normal in $\hat{G}$ one sees that $K$ is a proper subgroup of $G$; therefore $$[\hat{G}:K] > 3 \mbox{ showing that }\hat{G}/K \cong {\bf S}_3.$$Consequently, $K \leqslant G \cong \mathbb{Z}_p \times \mathbf{D}_4$ has order $4p$ and this implies that $$K = P \times V \mbox{ where } |V|=4 \mbox{ and } P=\langle a \rangle \cong \mathbb{Z}_p$$As $P$ is a normal Sylow $p$-subgroup of $K$ and $K \trianglelefteq \hat{G},$ we have that $P \trianglelefteq \hat{G}.$ Now, by a well-known result due to Zassenhaus, we conclude the existence of a subgroup $H$ of $\hat{G}$ of order 24 such that $\hat{G} = P \rtimes H.$  Observe that $$G \leqslant {\bf C}_{\hat{G}}(P) \unlhd \hat{G} \mbox{ and therefore } \mathbf{C}_{\hat{G}}(P) = \hat{G}.$$In other words, $P$ is a subgroup of the center of $\hat{G}$ and consequently $\hat{G} = P \times H.$ Finally, as $\hat{G}$ acts on $X$ with signature $(0; 2,3,4p),$ one has that $\hat{G}$ 
can be generated by an element of order $2$ and an element of order $3$. However, all the elements of $\hat{G}$ of these orders are contained in $H,$ a contradiction. This proves the claim.

\s

{\bf 2.} Assume $X \in \mathscr{M}_{2(p-1)}$ has a group of automorphisms $G$ of order $12p$. By proceeding as in the previous case, one obtains that the signature of the action of $G$ on $X$ has the form  \begin{equation}\label{gorro}(0; m_1, m_2, m_3) \mbox{ where }\tfrac{1}{m_1}+\tfrac{1}{m_2}+\tfrac{1}{m_3}=\tfrac{2}{3}+\tfrac{1}{2p}.\end{equation}

%
%
%
%

\s

{\it Claim.} The signature of the action of $G$ on $X$ is $(0; 2,6,2p)$ or $(0; 3,3,2p).$

\s

Let $v$ be the number of periods that are equal to two. Clearly,  $v \neq 2,3.$ Assume $v=1.$ Then, the signature of the action is $(0; 2, m_2, m_3)$ where $$\tfrac{1}{m_2}+\tfrac{1}{m_3}=\tfrac{1}{6}+\tfrac{1}{2p}.$$Without of generality, we can assume $m_2 \leqslant m_3.$ Note that $m_2 \neq 3.$ If $m_2=4$ then $$m_3=\tfrac{12p}{6-p} \mbox{ and therefore } p=5 \mbox{ and } m_3=60.$$Note that in such a case $G \cong \mathbb{Z}_{60}$ would act with signature $(0; 2,4,60);$ a contradiction. Similarly, if $m_2=5$ then $p=5$ and the signature of the action is $(0; 2,5,15).$ However, no group of order 60 acts with this signature (see, for instance, \cite{K}). Besides, if $m_2 > 6$ then $m_3 < 2p$ and hence $m_3 =12$ or $m_3=p.$ In the former case, one obtains that $m_2$ is not an integer, whereas in the latter one necessarily $p=5$ and the signature of the action is $(0; 2,10,15).$ But, no group of order 60 acts with this signature (see, for instance, \cite{K}).  All the above implies that the signature of the action is $(0; 2,6,2p).$

\s

Now, assume $v=0.$ Let $u$ be the number of periods that are equal to three. Clearly $u \neq 3.$ Assume that $u=1.$ Then the signature of the action is $$(0; 3, m_2, m_3) \mbox{ where }m_2, m_3 \geqslant 4 \mbox{ satisfy } \tfrac{1}{m_2}+\tfrac{1}{m_3}=\tfrac{1}{3}+\tfrac{1}{2p}.$$Observe that $m_2 < 2p$ and consequently $m_3 \in \{4,6,12,p\}.$ It can be seen that for each possible value of $m_3$ one obtains that $m_2$ is not an integer. Now, assume that $u=0.$ Note that if each period is greater or equal to 5 then \eqref{gorro} says that $p <0.$ Thus, we can assume $m_1=4.$ If $m_2=4$ then $m_3$ is not an integer, and we can assume $m_2, m_3 \geqslant 5.$ However, this contradicts \eqref{gorro}. All the above says that $u=2$ and the signature of the action is $(0; 3,3,2p).$ 

\s

{\it Claim.} If $G$ acts on $X$ with signature $(0; 2,6,2p)$ then $G \cong \mathbf{D}_3 \times \mathbf{D}_p.$

\s

As the claim holds for $p=5$ and $p=11$ (see, for instance,  \cite{K}), hereafter we assume that $p$ differs from these two values. Let $$\Delta_{(0; 2,6,2p)}=\langle \gamma_1, \gamma_2, \gamma_3 : \gamma_1^2=\gamma_2^6=\gamma_3^{2p}=\gamma_1\gamma_2\gamma_3=1\rangle$$and  assume that the action of $G$ on $X$ is represented by the ske $\theta: \Delta_{(0; 2,6,2p)} \to G.$ If $P \cong \mathbb{Z}_p$ is the normal $p$-Sylow subgroup of $G$ and $\pi : G \to G/P$ stands for the canonical projection, then  the composite map $$\eta:=\pi \circ \theta: \Delta_{(0; 2,6,2p)} \to G \to H:=G/P$$is a group epimorphism  onto a group $H$ of order 12. Set $$x=\eta(\gamma_1), \, y=\eta(\gamma_2) \mbox{ and } z=\eta(\gamma_3).$$We recall that $x^2=y^6=z^{2}=1$ and $xyz=1.$ Note that if $x=1$ then $y=z^{-1}$ and hence $y$ and $z$ are involutions, contradicting the surjectivity of $\eta.$ It follows that $x$ is an involution. Similarly, if $z=1$ then $y$ is an involution and $\eta$ is non-surjective. It follows that $z$ is an involution. As a consequence, being $H$ a group of order 12 generated by two involutions, we obtain that  $H \cong \mathbf{D}_6$. Now, the Schur-Zassenhaus theorem implies that  $$G \cong \langle a, r, s: a^p=r^6=s^2=(sr)^2=1, rar^{-1}=a^u, sas=a^{v} \rangle$$where $u^6\equiv 1 \mbox{ mod } p$ and $v=\pm 1.$ A straightforward computation shows that $u=\pm 1$ and $G$ is isomorphic to $\mathbb{Z}_p \times \mathbf{D}_6,  \mathbf{D}_{6p}$ or $\mathbf{D}_p \times \mathbf{D}_3.$ Clearly the former group does not act with signature $(0; 2,6,2p)$. The same conclusion holds for the second group, as its elements of order $6$ and $2p$ belong to its index two normal subgroup. This proves the claim.

\s

{\it Claim.} If $G$ acts on $X$ with signature $(0; 3,3,2p)$ then $p \equiv 1 \mbox{ mod }3$ and $G$ is isomorphic to the nontrivial semidirect product $\mathbb{Z}_p \rtimes \mathbf{A}_4.$

\s

As the claim holds for $p=5$ and $p=11$ (see, for instance, \cite{K}), we can assume $p$ different from these values. By proceeding as in the previous claim, the fact that $\mathbf{A}_4$ is the unique group of order 12 that can be generated by two elements of order 3 implies that $G/P \cong \mathbf{A}_4,$ where $P$ is the normal $p$-Sylow subgroup of $G.$ It follows that $G$ is isomorphic to  $\mathbb{Z}_p \times \mathbf{A}_4$ or to  \begin{multline}\label{oh}
\mathbb{Z}_p \rtimes_3 \mathbf{A}_4= \langle a, x,y,z: a^p=x^2=y^2=z^3= (xy)^2=1, \\zxz^{-1}=y, zyz^{-1}=xy, 
[a,x]=[a,y]=1, zaz^{-1}=a^r\rangle
\end{multline}where $r$ is a primitive third root of unity (and therefore $p \equiv 1 \mbox{ mod }3$). The claim follows after noting that the former group cannot be generated by elements of order 3.

\s

{\it Claim.} If $X$ is endowed with a group of automorphisms  $G$ isomorphic to $\mathbf{D}_p \times \mathbf{D}_3$, then $X \cong Y_p$ and  $\mbox{Aut}(X)=G$.

\s

We consider the group $G \cong \mathbf{D}_p \times \mathbf{D}_3$ with the presentation $$\langle  R, S, r, s: R^p=S^2=(SR)^2=r^3=s^2=(sr)^2=[R,r]=[R, s]=[S,r]=[S,s]=1 \rangle.$$By the previous claims, the signature of the action of $G$ on $X$ is $(0; 2,6,2p).$ Observe that the involutions of $G$ are $$sr^i, SR^j, sr^iSR^j \mbox{ for } 0 \leqslant i \leqslant 2 \mbox{ and }0 \leqslant j \leqslant p-1,$$whereas the elements of order six are $$r^kSR^l \mbox{ for } k=1,2 \mbox{ and } 0 \leqslant l \leqslant p-1.$$A routine computation shows that if the product of an element $g_1$  of order two and an element $g_2$  of order six has order $2p$ then $$g_1=sr^iSR^j \mbox{ and } g_2=r^kSR^l \mbox{ with } l \neq j.$$It follows that the action of $G$  on $X$  is represented by a ske  $$\theta_{i,j,k,l}=(sr^iSR^j, r^kSR^l, (sr^iSR^jr^kSR^l)^{-1})$$with $i,j,k,l$ as before and $l \neq j$. After considering the automorphism of $G$ given by $$r \mapsto r, s \mapsto sr^{-i}, R  \mapsto R, S \mapsto SR^{-j}$$we see that $\theta_{i,j,k,l}$ is $G$-equivalent to the ske $$(sS, r^kSR^{l-j}, (sSr^kSR^{l-j})^{-1}) \sim (sS, rSR,srR^{-1}).$$Consequently, $X$ is uniquely determined, up to isomorphism. As  $$(x,y) \mapsto  (\omega_3 x, y), \,\, (x,y) \mapsto (\tfrac{1}{x}, -\tfrac{y}{x^3}), \,\, (x,y) \mapsto (x, \omega_p y) \,\,\mbox{ and } \,\, (x,y)\mapsto (-x, \tfrac{1-x^6}{y})$$generate a group of automorphisms of $Y_p$ isomorphic to $\mathbf{D}_3 \times \mathbf{D}_p$, the uniqueness implies that $X \cong Y_p$. Finally, as the signature $(0; 2,6,2p)$ is maximal \cite{Sing72} we conclude that the full automorphism group of $X$ is $G.$

\s

{\it Claim.} If $p \equiv 1 \mbox{ mod }3$ and $X$ is endowed with a group of automorphisms $G$ isomorphic to $\mathbb{Z}_p \rtimes_3 \mathbf{A}_4$, then $X \cong Z_p$ and $\mbox{Aut}(X)=G.$

\s
We consider $G \cong \mathbb{Z}_p \rtimes_3 \mathbf{A}_4$ with the presentation \eqref{oh}, and recall that the signature of the action of $G$ on $X$ is $(0; 3,3,2p).$ Observe that the elements of $G$  of order $2p$ are $a^nx, a^ny, a^nxy$ for $1 \leqslant n \leqslant p-1,$ and the ones of order 3 are $$a^ix^jy^kz^s \mbox{ where }  i \in \{0, \ldots, p-1\}, j,k \in \{0, 1\}  \mbox{ and } s=\pm 1.$$

Let $\theta : \Delta_{(0; 3,3,2p)} \to G$ such that $\theta =(g_1, g_2, g_3)$ be a ske representing the action of $G$ on $X$.  After applying a suitable automorphism of $G$, we may assume that $g_3=ax$ and therefore $\theta$ is $G$-equivalent to either $$\theta_1=(g_1, a^ix^jy^kz, ax) \, \, \mbox{ or } \,\, \theta_2=(g_1, a^ix^jy^kz^2, ax)$$with $i,j,k$ as before. Now, the equality $$a^nx^my^l \cdot a^ix^jy^kz \cdot (a^nx^my^l)^{-1}=z \mbox{ where }n(r-1)=i, m=k-j \mbox{ and } l=k$$together with the fact that $[a^nx^my^l, ax]=1$ show that $\theta_1$ is $G$-equivalent to  $(a^{-1}xz^2,z,ax).$ Analogously, one sees that $\theta_2$ is $G$-equivalent to $(a^{-1}xz,z^2,ax).$ It then follows that, if $$S_i:=\mathbb{H}/\mbox{ker}(\theta_ i) \mbox{ for } i = 1,2$$ is the Riemann surface determined by $\theta_i$ then $X$ is isomorphic to $S_1$ or $S_2.$ As a consequence of the results of \cite{Sing72}, one has that $$\mbox{if } \mbox{Aut}(S_i) \neq G \mbox{ then } [\mbox{Aut}(S_i):G]=2.$$We claim that this situation does not occur. Indeed, the fact that each automorphism of $G$ must send $a$ to a power of it coupled with the relation $zaz^2=a^r$ show that there is no automorphism $\Phi$ of $G$ satisfying $$\Phi(z)=a^{-1}xz^2 \mbox{ nor } \Phi(z^2)=a^{-1}xz.$$ Now, following the results of \cite[Case N8]{Bu03}, we conclude that the action of $G$ on  $S_i$ does not extend, and hence the full automorphism group of $S_i$ agrees with $G,$ as claimed. Finally, if  $N(\Delta_{(0; 3,3,2p)})$ stands for the normaliser of $\Delta_{(0; 3,3,2p)}$ in $\mbox{Aut}(\mathbb{H}),$ then  $$N(\Delta_{(0; 3,3,2p)})/\Delta_{(0; 3,3,2p)}  \cong \mathbb{Z}_2$$acts faithfully by conjugation on  $\{\theta_1, \theta_2\}$ by identifying Riemann surfaces that are isomorphic. Thereby $S_1 \cong S_2$ and $X$ is uniquely determined, up to isomorphism. 

\s

Let $\lambda$ and $\mu$ be the integers satisfying $1+r+r^2=\lambda p$ and $r^3=1+\mu p$. A routine computation shows that transformations $$(x,y) \mapsto (x, \omega_py), \, (x,y) \mapsto (-x,y), \, (x,y) \mapsto (-\tfrac{1}{x}, \tfrac{y}{x^{4\lambda}})$$ and $$(x,y) \mapsto (i\tfrac{1-x}{1+x},  \tfrac{cy^{r^2}}{(x^2+1)^{2r \mu}(x^2-1)^{2\mu}(x+1)^{4\lambda}}) \mbox{ where } c^p=-4^{r+2r^2},$$  generate a group of automorphisms of $Z_p$ isomorphic to $\mathbb{Z}_p \rtimes_3 \mathbf{A}_4$. Therefore, the uniqueness of $X$ implies that $X \cong Z_p$, as claimed. This ends the proof of the theorem.

\subsection*{Proof of Theorem \ref{theobig}}
Let $\lambda \geqslant 4$ be an integer and assume that $X \in \mathscr{M}_{2(p-1)}$ has a group of automorphisms $G$ of order $4 \lambda p.$  Following \cite[p. 77]{accola}, the signature of the action of $G$ is $$\begin{array}{lcl}
(1) \,\, (0; 2,2,2,3), & \,\,\,\,\,  & (5) \,\, (0; 2,6,k) \mbox{ where } 6 \leqslant k \leqslant 11, \\
(2) \,\, (0; 2,3,k) \mbox{ where } k \geqslant 7, & \,\,\,\,\, & (6)\,\,(0; 2,7,k) \mbox{ where } 7 \leqslant k \leqslant 9, \\
(3) \,\, (0; 2,4,k) \mbox{ where } k \geqslant 5, & \,\,\,\,\, & (7) \,\, (0; 3,3,k) \mbox{ where } 4 \leqslant k \leqslant 11, \mbox{ or}\\
(4) \,\, (0; 2,5,k) \mbox{ where } 5 \leqslant k \leqslant 19, & \,\,\,\,\, & (8) \,\, (0; 3,4,k)\mbox{ where } 4 \leqslant k \leqslant 5.
\end{array}$$

As the theorem holds for $p=5$ and $p=7$ (see \cite{Conder}), hereafter we assume $p \geqslant 11.$

\s

{\bf 1.} $G$ does not act  with signature $(0; 2,2,3,3)$ since otherwise the Riemann-Hurwitz formula says that $\lambda=3-\frac{18}{4p} < 3.$ 

\s

{\bf 2.} If $G$ acts with signature $(0; 2,3,k)$ then $3$ divides $4p\lambda$ and therefore $\lambda=3 \lambda'$ for some $\lambda' \geqslant 2.$ The Riemann-Hurwitz formula reads $$4p-6=12p\lambda'(\tfrac{1}{6}-\tfrac{1}{k}) \, \iff \, k=\tfrac{6p\lambda'}{(\lambda'-2)p+3}.$$Note that if $\lambda'=3$ then $k$ is not an integer. Besides, if $\lambda' \geqslant 4$ then $$\tfrac{6\lambda'}{(\lambda'-2)p+3} \in \mathbb{Z} \, \implies p \leqslant 3(2+\tfrac{4}{\lambda'-2}) \leqslant \tfrac{21}{2};$$ a contradiction. It follows that $\lambda'=2$ and $G$ has order $24p$. However, this situation does not occur. Indeed, if $p \neq 11, 23$ then the $p$-Sylow subgroup of $G$ and a $2$-Sylow subgroup of $G$  generate a subgroup of order $8p$. But, by Theorem \ref{teo8p12p}, there is no compact Riemann surface of genus $g_p$ endowed with a group of automorphisms of order $8p$ with full automorphism group of order $24p.$ The conclusion also holds for $p=11, 23$ as can be seen from the database \cite{Conder}.

\s

{\bf 3.} If $G$ acts with signature $(0; 2,4,k)$ then the Riemann-Hurwitz formula reads $$4p-6=4p\lambda(\tfrac{1}{4}-\tfrac{1}{k}) \, \iff \, k=\tfrac{4p\lambda}{p(\lambda-4)+6}.$$Note that $\lambda \geqslant 5.$ It follows that $\tfrac{4\lambda}{p(\lambda-4)+6}$ is an integer and therefore  $p \leqslant 4+\tfrac{10}{\lambda-4}$. We see that necessarily $\lambda=5$ and $p=11,13.$ However, in such cases $k \notin \mathbb{Z}.$ 

\s

For the sake of conciseness, we omit the remaining cases as they are discarded analogously as in cases {\bf 2} and {\bf 3}.

\subsection*{Proof of Proposition \ref{lema1}}

Assume $X \in \mathscr{M}_{2(p-1)}$ has a group of automorphisms $G$ of order $4p$. Let $(h; m_1, \ldots, m_r)$ be the signature of the  corresponding action, and write $B=\Sigma_{i=1}^r(1-1/m_i).$

\s

If $h \geqslant 1$ then the Riemann-Hurwitz formula implies $4p-6 \geqslant 4pB,$ showing that $B =0.$ If $h=1$ then $2p=3,$ whereas if $h \geqslant 2$ one obtains that $p < 0.$ Thus, the signature of the action of $G$ on $X$ is $(0; m_1, \ldots, m_r)$ and $B=3-\tfrac{3}{2p}.$ 

\s

The fact that $B \leqslant \tfrac{r}{2}$ shows that $r$ is equal to $3,4$ or $5.$  Since $p \geqslant 5$ is prime, one can see easily that the latter case is impossible.

\s

We assume $r=4.$ Then \begin{equation}\label{ec1}\tfrac{1}{m_1}+\tfrac{1}{m_2}+\tfrac{1}{m_3}+\tfrac{1}{m_4}=1+\tfrac{3}{2p}.\end{equation}Let $v$ be the number of periods that are equal to 2. Clearly, $v \neq 0, 4.$ If $v=1$ then $$\tfrac{1}{m_2}+\tfrac{1}{m_3}+\tfrac{1}{m_4}=\tfrac{1}{2}+\tfrac{3}{2p}$$ contradicting the fact that $m_2, m_3, m_4 \geqslant 4.$ Similarly, if $v=3$ then \eqref{ec1} implies that $m_4$ is negative. All the above implies that the signature of the action is $$(0; 2,2, m_3, m_4) \mbox{ where } m_3, m_4 \in \{4,p,2p,4p\} \mbox{ and } \tfrac{1}{m_3}+\tfrac{1}{m_4}=\tfrac{3}{2p}.$$ The solutions of the equation above are $$m_3=p, m_4=2p \mbox{ for } p \geqslant 5 \mbox{ and } m_3=4, m_4=20 \mbox{ for }p=5.$$

We observe that the sporadic solution for $p=5$ is not realised. Indeed, if there is an action of a group of order 20 with signature $(0; 2,2,4,20)$ then $G$ is cyclic, and there are two involutions and an element of order 20  whose product has order 4; a contradiction.

\s

We now assume $r=3.$ Then \begin{equation*}\label{ec2}\tfrac{1}{m_1}+\tfrac{1}{m_2}+\tfrac{1}{m_3}=\tfrac{3}{2p}.\end{equation*}By proceeding as in the previous case, one sees that the number of periods that are equal to $2p$ is equal to $0$ or $3$. The latter case yields the solution $(0; 2p, 2p, 2p).$ In the former case, we note that at least one period must be equal to $4p,$ and thus $$\tfrac{1}{m_1}+\tfrac{1}{m_2}=\tfrac{5}{4p}$$The unique solution of this equation is $m_1=p$ and $m_2=4p$.

\s

Summing up, the possible signatures are $$s_1:=(0; 2,2,p,2p), \, s_2:=(0; 2p, 2p, 2p) \, \mbox{ and }\, s_3:=(0; p, 4p, 4p).$$ 

If the latter signature is realised then $G$ is cyclic. Conversely, if we write $\langle a \rangle \times \langle b \rangle \cong \mathbb{Z}_{p} \times \mathbb{Z}_4$ then the generating vector $(a, ab, a^{-2}b^{-1})$ shows that this group  acts with signature $s_3$.

\s

By the Sylow's theorems, if $p \equiv 3 \mbox{ mod } 4$ then $G$ is isomorphic to either $$\mathbb{Z}_{4} \times \mathbb{Z}_p, \,\, \mathbb{Z}_p \times \mathbb{Z}_2^2, \,\, \mathbf{D}_{2p} \,\, \mbox{ or }\,\, \langle a,b : a^p=b^4=1, bab^{-1}=a^{-1}\rangle=\mathbb{Z}_p \rtimes_2 \mathbb{Z}_4,$$and if $p \equiv 1 \mbox{ mod 4}$ then, in adition to these groups, $G$ can be isomorphic to $$\langle a,b : a^p=b^4=1, bab^{-1}=a^{r}\rangle=\mathbb{Z}_p \rtimes_4 \mathbb{Z}_4$$where $r$ is a primitive $4$-th root of unity in the field of $p$ elements. 

\s

Since the product of three elements of order $2p$ of $\mathbb{Z}_{4} \times \mathbb{Z}_p$ and $\mathbb{Z}_p \rtimes_2 \mathbb{Z}_4$ cannot be equal to the identity, we see that these groups do not act with signature $s_2$. The same conclusion is obtained for $\mathbf{D}_{2p}$ and $\mathbb{Z}_p \rtimes_4 \mathbb{Z}_4$, after noticing that elements of order $2p$ cannot generate the former, whereas the latter does not have elements of such an order. It follows that if $G$ acts with signature $s_2$ then $$G \cong \mathbb{Z}_p \times \mathbb{Z}_2^2=\langle t, x,y : t^p=x^2=y^2=(xy)^2=[t,x]=[t,y]=1\rangle.$$Conversely, the tuple $(tx, ty, t^{-2}xy)$ shows that this group acts with signature $s_2$.

\s

Finally, the facts that $\mathbb{Z}_p \times \mathbb{Z}_4$ and $\mathbb{Z}_p \rtimes_2 \mathbb{Z}_4$ have a unique involution each, and that $\mathbb{Z}_p \rtimes_4 \mathbb{Z}_4$  does not have elements of order $2p$ show that these groups do not act with signature $s_1.$ Conversely, if we write $\mathbf{D}_{2p}=\langle r, s : r^{2p}=s^2=(sr)^2=1\rangle$, then the generating vectors  $$(x,y,t,xyt^{-1}) \mbox{ and } (s,sr^{2p-3},r^2,r)$$show that $\mathbb{Z}_p \times \mathbb{Z}_2^2$ and $\mathbf{D}_{2p}$ act with signature $s_1.$

\subsection*{Proof of Theorem  \ref{teoA}} Let $p \geqslant 5$ be a prime number, and consider $$\Delta_{(0; 2,2,p,2p)}=\langle \gamma_1, \gamma_2, \gamma_3, \gamma_4 : \gamma_1^2= \gamma_2^2= \gamma_3^p= \gamma_4^{2p}=\gamma_1 \gamma_2 \gamma_3 \gamma_4=1 \rangle.$$

\subsection*{1.} Consider the group $$G=\langle t, x,y : t^p=x^2=y^2=(xy)^2=[t,x]=[t,y]=1 \rangle \cong \mathbb{Z}_p \times \mathbb{Z}_2^2.$$The existence of the family $\mathcal{F}_1$ is guaranteed by the fact that $$\theta_0: \Delta_{(0; 2,2,p,2p)} \to G \mbox{ given by }\theta_0=(x,y,t,xyt^{-1})$$is a ske. Following \cite{B90}, the irreducible components of $\mathcal{F}_1$ are in bijective correspondence with the topological classes of actions of $G$ with signature $s_1$. Thus, we need to show that there is only one among such classes. Assume that the ske $\theta=(g_1, g_2, g_3, g_4)$ represents an action of $G$ in genus $g_p$ with signature $s_1.$ Observe that  $g_1 \neq g_2$ and then, up to an isomorphism of $G$, we can assume that $g_1=x$ and $g_2=y.$ Besides, after considering an appropriate automorphism of $G$, we can choose $g_3$ to be equal to $t$ and therefore $\theta$ is equivalent to $\theta_0,$ as desired. 

\s

We now proceed to provide an algebraic description of each $X \in \mathcal{F}_1.$
The regular covering map given by the action of $\langle t \rangle$ on $X$ $$X \to Y=X/\langle t \rangle \cong \mathbb{P}^1$$ramifies over six values. As $G/\langle t \rangle\cong \mathbb{Z}_2^2$ acts on $Y$ with signature $(0; 2,2,2)$, we can assume that such branch values are $0, \infty, \pm 1, \pm \mu,$ for some $\mu \in \mathbb{C}$ such that $\mu \neq 0, \pm 1.$ It follows that $X$ is  represented by the algebraic curve $$y^p=x^{\epsilon_0}(x-1)^{\epsilon_1}(x+1)^{\epsilon_2}(x-\mu)^{\epsilon_3}(x+\mu)^{\epsilon_4}$$for some integers $1 \leqslant \epsilon_0, \ldots, \epsilon_4 < p$ such that $\epsilon_0+\cdots+\epsilon_4$ is coprime with $p.$

We denote by $\pi : X \to X/G$ the regular covering map given by the action of $G$ on $X$. As noticed above, this action may be represented by the ske $$(x,y,t^{-2},xyt^{2}).$$

Let $Q$ be the branch value of $\pi$ marked with $p$. The fibre $\pi^{-1}(Q)$ consists of four points $P_1, P_2,P_3, P_4$ with $G$-stabiliser $\langle t \rangle$. We may assume that (see, for instance, \cite[Lemma 3.3]{B22})  $$\mbox{rot}(P_1, t^{-2})=\omega_{p} \mbox{ showing that }\epsilon_1=p-2.$$ The fact $G$ is abelian implies that $\epsilon_1=\epsilon_2=\epsilon_3=\epsilon_4.$ By proceeding analogously, one sees that $\epsilon_0=4.$  Hence, $X$ is isomorphic to the normalisation of \begin{equation}\label{lapizqq}y^p=x^{4}[(x^2-1)(x^2-\mu^2)]^{p-2}\end{equation}(see, for instance, \cite[Proposition 2.45]{libroge}). Observe that the transformations $$\mathbf{a}(x,y)=(x,\omega_p y), \,\, \mathbf{x}(x,y)=(-x,y)\, \mbox{ and } \,\mathbf{y}(x,y)=(\tfrac{\sqrt{t}}{x}, \tfrac{ty}{x^4})$$are automorphisms of \eqref{lapizqq} and $\langle \mathbf{a} \rangle \times \langle \mathbf{x}, \mathbf{y}\rangle \cong \mathbb{Z}_p \times \mathbb{Z}_2^2.$

\s

By the results of \cite{Sing72}, no Fuchsian group of dimension one contains $\Delta_{(0; 2,2,p,2p)}$. Hence, up to finitely many exceptional cases, the full automorphism group of each member of $\mathcal{F}_1$ is $G$. 
\s

Finally, note that by taking $\mu=\sqrt{-1}$ one sees that $X_p$ belongs to $\mathcal{F}_1.$

\subsection*{2.}Consider the group $$G=\langle r, s : r^{2p}=s^2=(sr)^2=1 \rangle \cong \mathbf{D}_{2p} .$$The existence of the family $\mathcal{F}_2$ is guaranteed by the fact that $$\theta_0: \Delta_{(0; 2,2,p,2p)} \to G \mbox{ given by }\theta_0=(s,sr^{2p-3},r^2,r)$$is a ske. Now, assume that the ske $ \theta=(g_1, g_2, g_3, g_4)$ represents an action of $G$ in genus $g_p$ with signature $s_1.$ Clearly,  $g_3=r^{2k}$ for some $k \in \{1, \ldots, p-1\}.$ Observe that, up to an automorphism of $G$ the form $r \mapsto r^u$, we can assume that $g_4=r.$ In addition, the surjectivity of $\theta$ implies that $g_1$ and $g_2$ are non-central, and therefore, after considering an automorphism of $G$ of the form $s \mapsto sr^{u},$ we may assume $g_1=s$. It follows that the action represented by $\theta$ is topologically equivalent to the action represented by $$\theta_k=(s,sr^{2p-2k-1}, r^{2k}, r) \mbox{ for some } k \in \{1, \ldots, p-1\}.$$  We recall that (see, for instance, \cite{B91}) the transformation $$\Phi : (\gamma_1, \gamma_2, \gamma_3, \gamma_4) \mapsto (\gamma_1, \gamma_3^{-1}, \gamma_3\gamma_2\gamma_3, \gamma_4)$$identifies skes representing actions that are topologically equivalent. The fact that $$(\Phi \circ \Phi) \cdot \theta_k = \theta_{2p-k},$$ shows that $\theta$ is topologically equivalent to $\theta_k$ for some  $k \in \{1, \ldots, \tfrac{p-1}{2}\}.$ In other words, the family $\mathcal{F}_2$ consists of at most $\tfrac{p-1}{2}$ irreducible components.

\s

We now proceed to provide an algebraic description of each $X \in \mathcal{F}_2.$
The regular covering map $X \to Y=X/\langle r^2 \rangle \cong \mathbb{P}^1$ given by the action of $\langle r^2 \rangle$ on $X$ ramifies over six values. As such branch values are invariant under the action of $G/\langle r^2 \rangle \cong \mathbb{Z}_2^2$, they can be assumed to be $\pm 1, \pm t, \pm \tfrac{1}{t}$ for some $t \in \mathbb{C}-\{0, \pm 1, \pm i\}$ (here we choose $\pm1$ as the ramification values associated to the branch points with $G$-stabiliser of order $2p$). 

It follows that $X$ is  represented by the algebraic curve $$y^p=(x-1)^{\epsilon_1}(x+1)^{\epsilon_2}(x-t)^{\epsilon_3}(x+t)^{\epsilon_4}(x-\tfrac{1}{t})^{\epsilon_5}(x+\tfrac{1}{t})^{\epsilon_6}$$for some integers $1 \leqslant \epsilon_1, \ldots, \epsilon_6 < p$ such that $\epsilon_1+\cdots+\epsilon_6$ is divisible by $p.$

We denote by $\pi : X \to X/G$ the regular covering map given by the action of $G$ on $X$, and assume that $X$ belongs to the component determined by  $$\theta_k=(s,sr^{p-2k-1},r^{2k},r)\mbox{ where } k \in \{1, \ldots, \tfrac{p-1}{2}\}.$$

Let $Q$ be the branch value of $\pi$ marked with $p$. Then, the fibre of $Q$ is formed by $P, sP, rP$ and $srP.$ Without loss of generality, one has that $$\mbox{rot}(P, r^{2k})=\omega_{p} \mbox{ showing that }\epsilon_3=k,$$and consequently  $\epsilon_4=p-k.$ By proceeding analogously, we obtain that $\epsilon_5=\epsilon_3, \epsilon_6=\epsilon_4, \epsilon_1=1$ and $\epsilon_2=p-1.$ Hence, $X$ is isomorphic to the normalisation of  \begin{equation}\label{lapiz12}y^p=(x-1)(x-t)^k(x-\tfrac{1}{t})^k(x+1)^{p-1}(x+t)^{p-k}(x+\tfrac{1}{t})^{p-k}\end{equation}

Observe that the transformations $$\mathbf{u}(x,y)=(x,\omega_p y) \, \mbox{ and } \,\mathbf{v}(x,y)=(\tfrac{1}{x}, -\tfrac{y}{x^3})$$are automorphisms of \eqref{lapiz12} of order $p$ and 2 that commute. Also, the  map $$\mathbf{s}(x,y) \mapsto (-x,\tfrac{1}{y}(x^2-1)(x^2-t^2)(x^2-\tfrac{1}{t^2}))$$is an automorphism of \eqref{lapiz12} of order two. Note that $\langle \mathbf{u}\mathbf{v}, \mathbf{s}\rangle \cong \mathbf{D}_{2p}$ and therefore these maps generate $G$. As argued in the proof of the previous proposition, the maximality of $\Delta_{(0; 2,2,p,2p)}$ shows that,  up to some finitely many exceptional cases, the full automorphism group of each member $X$ of $\mathcal{F}_2$ is $G.$

\s

Finally, note that by taking $t=\omega_3$ and $k=1$ one sees that $Y_p$ belongs to $\mathcal{F}_2$.

\subsection*{Proof of Corollary \ref{coro1}} If $\mathcal{S} \in \mathcal{F}_1 \cap \mathcal{F}_2$ then the automorphism group of $\mathcal{S}$ has order strictly than $4p$ (and hence, by Theorem \ref{teo8p12p} and Theorem \ref{theobig}, it is isomorphic to either $X_p, Y_p$ or $Z_p$) and has subgroups isomorphic to $\mathbb{Z}_p \times \mathbb{Z}_2^2$ and $\mathbf{D}_{2p}.$ However, none of such Riemann surfaces satisfies this property. Similarly, if $\mathcal{S} \in \partial(\mathcal{F}_1)$ is nonisomorphic to $X_p$ then it is isomorphic to $Y_p$ or $Z_p.$ But the automorphism groups of these surfaces do not have subgroups isomorphic to $\mathbb{Z}_p \times \mathbb{Z}_2^2.$ We argue analogosuly for $\mathcal{F}_2.$

\subsection*{Proof of Theorem \ref{teoD}} 

Assume $X \in \mathscr{M}_{2(p-1)}$ has a group of automorphisms $G$ isomorphic to
$$\mathbb{Z}_p \times \mathbb{Z}_2^2=\langle t, x,y : t^p=x^2=y^2=(xy)^2=[t,x]=[t,y]=1\rangle$$acting with signature $s_2.$ The elements of order $2p$ are $t^ix, t^iy, t^ixy$ where $1 \leqslant i \leqslant p-1.$ The action of $\mbox{Aut}(G)\cong \mathbb{Z}_{p-1} \times \mathbf{S}_3$ on the skes of $G$ of signature $s_2$ has $p-2$ orbits represented by $$\theta_j: \Delta_{(0; 2p, 2p, 2p)} \to G \mbox{ given by }\theta_j=(tx, t^jy, t^{-1-j}xy) \mbox{ where } j \in \{1, \ldots, p-2\}.$$Moreover, the action of the braid group $\mathscr{B}$ on $\{\theta_1, \ldots, \theta_{p-2}\}$ show that the topological classes of actions are determined by $$\theta_j \sim \theta_{\frac{1}{j}} \sim \theta_{-1-j} \sim \theta_{\frac{-j}{1+j}} \sim \theta_{-1-\frac{1}{j}} \sim \theta_{\frac{-1}{1+j}}$$(see, for instance, \cite[\S4.1.1]{B91} for generators of $\mathscr{B}$). In other words, $\theta_{j_1}$ is topologically equivalent to $\theta_{j_1}$ if and only if  $j_1=\sigma(j_2)$ where $\sigma \in \langle z \mapsto 1/z, z \mapsto -1/(1+z)\rangle,$ with the subindices taken modulo $p.$ Let $R_j$ be the Riemann surface $\mathbb{H}/\mbox{ker}(\theta_j)$. 

\s

{\it Claim.}  $R_1 \cong R_{p-2} \cong R_{\frac{p-1}{2}} \cong X_p$.

\s

We recall that, as shown in the proof of Theorem \ref{teo8p12p}, \begin{equation*}\label{g1}\mbox{Aut}(X_p) \cong \langle a,r,s : a^p=r^4=s^2=(sr)^2=[a,r]=[a,s]=1\rangle \cong \mathbb{Z}_p \times \mathbf{D}_4,\end{equation*}and its action on $X_p$  is represented by the ske $$\Theta : \Delta_{(0; 2, 2p, 4p)} \to \mathbb{Z}_p \times \mathbf{D}_4 \mbox{ given by } (\gamma_1, \gamma_2, \gamma_3) \mapsto (s, a^{-1}sr, ar^{-1}).$$Note that  $\Gamma:=\langle \gamma_2, \gamma_3\gamma_2\gamma_3^{-1},\gamma_3^2 \rangle \cong \Delta_{(0; 2p, 2p, 2p)}$ and the restriction of $\Theta$ to $\Gamma$ is given by $$(\gamma_2, \gamma_3\gamma_2\gamma_3^{-1},\gamma_3^2) \mapsto (a^{-1}sr, a^{-1}sr^3, a^2r^2).$$If we let $t':=a^{-1}, x':=sr$ and $y':=sr^3$ then the restriction above yields an action of $$H:=\langle t', x', y'\rangle \cong \mathbb{Z}_p\times \mathbb{Z}_2^2$$that is $H$-equivalent to the one of $\theta_1.$ Hence, $X_p$ is isomorphic to $R_1.$ 
As $\theta_1, \theta_{p-2}$ and $\theta_{\frac{p-1}{2}}$ are topologically equivalent, we obtain that the action of $G$ on $R_{p-2}$ and on $R_{\frac{p-1}{2}}$ also extend to an action of $\mathbb{Z}_p \times \mathbf{D}_4$, and the claim follows from the uniqueness of $X_p.$

\s

{\it Claim.} If $p \equiv 1 \mbox{ mod }3$ and $2 \leqslant r \leqslant p-2$ satisfies $r^3 \equiv 1 \mbox{ mod }p,$ then $R_r \cong R_{r^2} \cong Z_p.$

\s

The automorphism group of $Z_p$ is isomorphic to   \begin{multline*}
\mathbb{Z}_p \rtimes_3 \mathbf{A}_4= \langle a, x,y,z: a^p=x^2=y^2=z^3= (xy)^2=1, \\zxz^{-1}=y, zyz^{-1}=xy, 
[a,x]=[a,y]=1, zaz^{-1}=a^r\rangle
\end{multline*}and its action on $Z_p$  is represented by the ske $$\Theta : \Delta_{(0; 3,3,2p)} \to \mathbb{Z}_p \rtimes_3 \mathbf{A}_4\mbox{ given by } (y_1, y_2, y_3) \mapsto (a^{-1}xz^2, z,az).$$

The restriction of $\Theta$ to $\Gamma':=\langle \gamma_2\gamma_3\gamma_2^{-1},\gamma_2^{-1}\gamma_3\gamma_2, \gamma_3\rangle \cong \Delta_{(0; 2p, 2p, 2p)}$ is given by $$(\gamma_2\gamma_3\gamma_2^{-1},\gamma_2^{-1}\gamma_3\gamma_2, \gamma_3) \mapsto (a^ry, a^{r^2}xy, ax).$$If we write $t':=a^{r}, x':=y$ and $y':=xy$ then the restriction yields an action of $$H':=\langle a', x', y'\rangle \cong \mathbb{Z}_p\times \mathbb{Z}_2^2$$that is $H$-equivalent to the one of $\theta_r.$ We conclude as in the previous claim. \s

It follows from the claims above if $1 \leqslant j \leqslant p-2$ is different from $1, p-2, \tfrac{p-1}{2}$ and distinct from the primitive third roots of unity in the case that $p \equiv 1 \mbox{ mod }3,$ then the full automorphism group of $R_j$ is $G$. Let $\mathcal{P}$ be this collection of subindices. 

\s

By the results of \cite{Sing72}, the normaliser of $\Delta_{(0; 2p, 2p, 2p)}$ is $\Delta_{(0; 2,3,4p)}$ and the action of $$\Delta_{(0; 2,3,4p)}/\Delta_{(0; 2p, 2p, 2p)} \cong \mathbf{S}_3$$by conjugation in $\{\theta_j : j \in \mathcal{P}\}$ identifies skes representing Riemann surfaces that are isomorphic. It follows that if $p \equiv 1 \mbox{ mod } 3$ then there are $|\mathcal{P}|/6+2$ pairwise nonisomorphic Riemann surfaces in $\mathcal{C}=\{R_1, \ldots, R_{p-2}\}$, and if $p \equiv 2 \mbox{ mod } 3$ then there are $|\mathcal{P}|/6+1$ pairwise nonisomorphic Riemann surfaces in $\mathcal{C}$. Summing up, there are $$  \tfrac{p+5}{6} \mbox{ if } p \equiv 1 \mbox{ mod } 3 \,\, \mbox{ and } \,\,
 \tfrac{p+1}{6} \mbox{ if }p \equiv 2 \mbox{ mod } 3$$pairwise nonisomorphic Riemann surfaces in $\mathscr{M}_{2(p-1)}$ with a group of automorphisms of order $4p$ acting with signature $s_2.$ 

\s

To finish the proof, observe that the action of $G$ on $R_j$ is also represented by the ske $(t^{-1}x, t^{-j}y, t^{1+j}xy)$. Now, we proceed as in the proof of the previous theorem, to see that the rotation number of $t$ at its fixed points with $G$-stabiliser $\langle t^{-1}x \rangle$ is $p-2$, at the ones with $G$-stabiliser $\langle t^{-j}x \rangle$ is $p-2j$, and at the ones $G$-stabiliser $\langle t^{1+j}xy \rangle$ are $2j+2$. Thus, $R_j$ is represented by the curve $$y^p=x^{2j+2}(x^2-1)^{p-2}(x^2+1)^{p-2j}$$and hence $R_j \cong C_j$ as desired.

\subsection*{Proof of Corollary \ref{normal}}

Following \cite[Section 2]{GG92} (see also \cite{HPRCR}), the family $\mathcal{F}_1$ is a non-normal subvariety of $\mathscr{M}_{2(p-1)}$ if and only if there exists $\mathcal{S} \in \mathcal{F}_1$ satisfying the following two properties.
\begin{enumerate}
\item  $\mbox{Aut}(\mathcal{S})$ has two non-conjugate subgroups $H_1$ and $H_2$ that are isomorphic to $\mathbb{Z}_p \times \mathbb{Z}_2^2.$
\item The actions of $H_1$ and $H_2$ on $\mathcal{S}$  are topologically equivalent.
\end{enumerate}

If $\mathcal{F}_1$ is non-normal, then by Corollary \ref{coro1}, necessarily $\mathcal{S}=X_p$. By Theorem \ref{teo8p12p}, the automorphism group of $X_p$ is isomorphic to $\mathbb{Z}_p \times \mathbf{D}_4$ and  it has exactly two non-conjugate subgroups isomorphic to $\mathbb{Z}_p \times \mathbb{Z}_2^2.$ Clearly, one of these subgroups must act with signature $s_1,$ whereas the remaining subgroup must act with signature $s_2$ as shown in Theorem \ref{teoD}. It follows that the corresponding actions are topologically inequivalent. We then conclude that each point of the complex one-dimensional family $\mathcal{F}_1$ is normal, and hence smooth. The same argument applies for each irreducible component of $\mathcal{F}_2.$

\subsection*{Proof of Theorem \ref{theoC}}

Assume $X \in \mathscr{M}_{2(p-1)}$ has a group of automorphisms $G$ isomorphic to$$\langle a, b : a^p=b^4=[a,b]=1\rangle \cong \mathbb{Z}_{4p}.$$By Lemma \ref{lema1}, the signature of the action is $s_3,$ and after considering a suitable automorphism of $G$, we see that the action of $G$ is represented by the ske $$\theta_j: \Delta_{(0; p,4p,4p)} \to G \mbox{ such that } \theta_j=(a, a^jb,a^{-j-1}b^{3})$$for some $j \in \{1, \ldots, p-2\}.$ Therefore, if we write $$\mathcal{Y}_j:=\mathbb{H}/\mbox{ker}(\theta_j) \mbox{ then }X \cong \mathcal{Y}_j \mbox{ for some }j \in\{1, \ldots, p-2\}.$$

\s

{\it Claim.} $\mathcal{Y}_{\frac{p-1}{2}} \cong X_p$ and $\mathcal{Y}_j \cong \mathcal{Y}_{p-1-j}$ for each $j \in \{1, \ldots, \tfrac{p-3}{2}\}.$

\s

There are exactly two triangle Fuchsian groups that contain $\Delta_{(0; p,4p,4p)}$ properly$$\Delta_{(0; p,4p,4p)} \leqslant \Delta_{(0; 2,2p,4p)} \leqslant \Delta_{(0; 2,3,4p)}$$(see \cite{Sing72}). This implies that if the full automorphism group of $\mathcal{Y}_j$ is different from $G$ then one of the following statements holds.

\begin{enumerate}
\item $\mbox{Aut}(\mathcal{Y}_j)$ has order $8p$ and acts with signature $(0; 2,2p,4p).$
\item $\mbox{Aut}(\mathcal{Y}_j)$ has order $12p$ and acts with signature $(0; 2,3,4p).$
\end{enumerate}The latter case is not realised, as a consequence of Theorem \ref{teo8p12p}. We recall that, as noticed in the proof of Theorem \ref{teo8p12p}, the action of $H \cong \mathbb{Z}_p \times \mathbf{D}_4=\langle a \rangle \times \langle r,s \rangle$ on $X_p$ is determined by $$\Theta: \Delta_{(0; 2,2p,4p)} \to H \mbox{ given by } (\gamma_1, \gamma_2, \gamma_3) \mapsto (s,a^{-1}sr,ar^{-1}).$$ Observe that $\Gamma:=\langle \gamma_2^2, \gamma_2^{-1}\gamma_3\gamma_2, \gamma_3  \rangle$ is isomorphic to $\Delta_{(0; p,4p,4p)}$ and the restriction$$\Theta|_{\Delta_{(0; p,4p,4p)}} \mbox{ is given by } (\gamma_2^2, \gamma_2^{-1}\gamma_3\gamma_2, \gamma_3) \mapsto (a^{-2}, ar, ar^{-1}).$$After considering the automorphism of $\Theta(\Delta_{(0; p,4p,4p)}) =\langle ar \rangle \cong \mathbb{Z}_{4p}$ given by $a \mapsto a^{\frac{p-1}{2}}$ and $r \mapsto r,$ we are in position to conclude that $$\Theta|_{\Delta_{(0; p,4p,4p)}}=\theta_{\frac{p-1}{2}} \mbox{ and hence }X_p \cong \mathcal{Y}_{\frac{p-1}{2}}.$$ Besides, observe that if $j \neq  \frac{p-1}{2}$ then there is no automorphism $\Phi$ of $G$ that satisfies $$\Phi(a)=a \mbox{ and } \Phi(a^jb)=a^{-1-j}b^3.$$This implies, by \cite[Case N8]{Bu03}, that $\mbox{Aut}(\mathcal{Y}_j)=G$ if and only if $j \neq \tfrac{p-1}{2}.$

\s

The normaliser of $\Delta_{(0; p,4p,4p)}$ is $\Delta_{(0; 2,2p,4p)}$ and the quotient group $$\Delta_{(0; 2,2p,4p)}/\Delta_{(0; p,4p,4p)} = \langle X_2\Delta_{(0; p,4p,4p)} \rangle \cong \mathbb{Z}_2$$acts faithfully by conjugation on $\{\theta_1, \ldots, \theta_{p-2}\}-\{\theta_{\frac{p-1}{2}}\}$. Two skes belong to the same orbit if and only if the Riemann surfaces they represent are isomorphic. Thus, the claim follows after noticing that the orbits  are $$\{\theta_j, \theta_{p-1-j}\} \mbox{ for }j \in \{1, \ldots, \tfrac{p-3}{2}\}.$$

The regular covering map $\mathcal{Y}_j \to \mathbb{P}^1$ given by the action of $\langle a \rangle$ on $\mathcal{Y}_j$ ramifies over six values. As such branch values are invariant under the action of $G/\langle a \rangle \cong \mathbb{Z}_4=\langle \hat{b} \rangle$, they can be assumed to be $\pm 1, \pm i$ (forming an orbit of length 4) and $0, \infty$ (being the fixed points of $\tilde{b}$).  It follows that $\mathcal{Y}_j$ is  represented by  $$y^p=x^{\epsilon_0}(x-1)^{\epsilon_1}(x-i)^{\epsilon_2}(x+1)^{\epsilon_3}(x+i)^{\epsilon_4}$$for some integers $1 \leqslant \epsilon_0, \ldots, \epsilon_4 < p$ such that $\epsilon_0+\cdots+\epsilon_4$ is coprime with  $p.$

We denote by $\pi_j$ the regular covering map given by the action of $G$ on $\mathcal{Y}_j$, which is represented by the ske  $$\theta_j=(a,a^jb,a^{-j-1}b^3) \sim (a^s, ab, a^{-1-s}b^{3})=:\psi_s$$where $sj\equiv 1 \mbox{ mod }p.$ If $P$ is a branch point of $\pi_j$ with $G$-stabiliser $\langle a \rangle$ and $P'$ is the branch point of $\pi_j$ with $G$-stabiliser $\langle ab\rangle$ then, without loss of generality, we may assume $$\mbox{rot}(P, a^s)=\omega_{p} \mbox{ and }\mbox{rot}(P', a^4)=\omega_{p}.$$Hence, we proceed as done in the theorem above to see that  $\mathcal{Y}_j \cong S_j$ as desired.

\subsection*{Group algebra decomposition of Jacobian varieties}\mbox{}

Group actions on Jacobian varieties have been successfully employed to
decompose them isogenously. More precisely, following Lange and Recillas \cite{LR04}, a $G$-action on a compact Riemann surface 
$X$ induces the so-called isotypical decomposition of the Jacobian variety $JX$
$$ JX \sim A_1 \times \cdots \times A_t$$
in terms of abelian subvarieties $A_i$ of $JX$ that are pairwise non-$G$-isogenous. Moreover, each
$A_i$ decomposes further as $B_i^{n_i} $
for some abelian subvariety $B_i$ of $JX$. Hence
$$ JX \sim B_1^{n_1} \times \cdots \times B_t^{n_t}$$
which is called the group algebra decomposition of $JX$ with respect to $G$. As proved by
Rodríguez and the first author in \cite{CR06}, the integers $n_i$  can be determined explicitly from algebraic
data coming from the rational irreducible representations of $G.$ Further, the dimension
of the factors $B_i$ depends heavily on how the group acts. Rojas in \cite{R07} gave an explicit formula
to compute such dimensions.

\s

It was also shown in \cite{LR04} and \cite{CR06} that if $H \leqslant
N \leqslant G$ are subgroups of $G$ with intermediate covering $$F :
X_H:=X/H \to X_N:=X/N,$$then the
Jacobians $JX_H$ and $JX_N$, as well as the (generalized) Prym
variety $P(X_H/X_N)$, defined as the orthogonal complement of
$F^*(JX_N)$ in $JX_H$, admit similar isogeny decompositions. In
fact, there exist non-negative integers $h_j$ and $p_j$ such
that
$$
JX_H \sim B_1^{h_1} \times \cdots \times B_r^{h_t}  \,\mbox{ and } \, P(X_H/X_N) \sim B_1^{p_1} \times \cdots \times B_r^{p_t}.
$$

It was further shown in \cite{CLRR} and in \cite{e} that the
Hecke algebra $\mathbb{Q}[H\backslash G/H]$; namely the subalgebra of
$\mathbb{Q}[G]$ consisting of the $\mathbb{Q}$-valued functions
on $G$ that are constant on each double coset of $H$ in $G,$ acts naturally on
$JX_H$.

\subsection*{Proof of Theorem \ref{jacofam}}\mbox{} We recall that if $H$ is a group of automorphisms of $X$ then $X_H$ stands for the quotient $X/H.$

{\bf 1.} Consider the complex irreducible representations of \begin{equation*}\label{tabaco}G=\mathbb{Z}_p \times \mathbb{Z}_2^2=\langle t, x,y : t^p=x^2=y^2=(xy)^2=[t,x]=[t,y]=1\rangle\end{equation*}given by \begin{equation}\label{tabaco1}V_1 : (t,x,y) \mapsto (\omega_p, 1, -1), V_2 : (t,x,y) \mapsto (\omega_p, -1, 1),V_3 : (t,x,y) \mapsto (\omega_p, -1, -1).\end{equation}Note that these representations are pairwise Galois non-asociated. Let $$W_j:=\oplus_{\mu \in H}V_j^{\mu} \mbox{ where } H=\mbox{Gal}(\mathbb{Q}(\omega_p)/\mathbb{Q})$$ denote the rational irreducible representation of $G$ associated to $V_j$. As recalled in the previous subsection, the group algebra decomposition of $JX$ for each $X \in \mathcal{F}_1$ with respect to $G$ has the form
$$JX \sim B_{W_1}\times B_{W_2} \times B_{W_3} \times P$$where $B_{W_j}$ is the abelian subvariety of $JX$ associated to $W_j$ and $P$ is the abelian subvariety of $JX$ given by the product of the factors corresponding to the remaining rational irreducible representations of $G$ (see \cite{CR06} and \cite{LR04}). We recall that the action of $G$ on $X$ is given by the ske $(x,y,t,xyt^{-1})$. By applying  \cite[Theorem 5.12]{R07} we have that $$\mbox{dim}B_{W_1}=\dim B_{W_2}=\tfrac{p-1}{2} \mbox{ and } \dim B_{W_3}=p-1.$$Thus, the dimension of $P$ is zero and hence $$JX \sim B_{W_1}\times B_{W_2} \times B_{W_3}.$$ Now, by applying \cite[Proposition 5.2]{CR06}, one obtains that $$JX_{\langle x \rangle} \sim B_{W_1}, \,\, JX_{\langle y \rangle} \sim B_{W_2} \, \mbox{ and } \, JX_{\langle xy \rangle} \sim B_{W_3}$$and the proof is done after noting that $x$ and $y$ are precisely the involutions of $X$ with exactly $2p$ fixed points.

\s

{\bf 2.}  Consider the complex irreducible representations of degree two of $$G=\mathbf{D}_{2p}=\langle r, s : r^{2p}=s^2=(sr)^2=1\rangle$$given by  $$V_1 : r \mapsto (\omega_{2p}, \bar{\omega}_{2p}) \, s \mapsto  \left( \begin{smallmatrix}
0 & 1 \\
1 & 0
\end{smallmatrix} \right) \, \mbox{ and } \,  V_2 : r \mapsto (\omega_{p}, \bar{\omega}_{p}), s \mapsto  \left( \begin{smallmatrix}
0 & 1 \\
1 & 0
\end{smallmatrix} \right)$$ Both representations are defined over their character fields $\mathbb{Q}(\omega_p+\bar{\omega}_p)$, have Schur index 1, and are Galois non-asociated. Let $G_j$ be the Galois group of $V_j$ and let $$W_j:=\oplus_{\mu \in G_j}V_j^{\mu}$$ be the rational irreducible representation of $G$ associated to $V_j$. Then the group algebra composition of $JX$ for each $X \in \mathcal{F}_2$ has the form$$JX \sim B_{W_1}^2\times B_{W_2}^2  \times P$$for some abelian subvariety $P$ of $JX.$ As the action is given by the ske $$(s,sr^{2p-2k-1}, r^{2k}, r) \mbox{ for some }k \in \{1, \ldots, \tfrac{p-1}{2}\}$$ we see that \cite{R07} implies that, independently of the value of $k,$ $$\mbox{dim}B_{W_1}=\dim B_{W_2}=\tfrac{p-1}{2}$$Thus, the dimension of $P$ is zero and hence $$JX \sim B_{W_1}^2\times B_{W_2}^2.$$ Let $\iota$ be a non-central involution of $X$. Then $\iota$ is conjugate to $s$ and, by applying the \cite[Proposition 5.2]{CR06}, one obtains that $$JX_{\langle \iota \rangle} \cong JX_{\langle s \rangle} \sim B_{W_1} \times  B_{W_2}.$$Note that $JX$ decomposes a product of four  subvarieties of dimension $\tfrac{p-1}{2}.$

\s

{\bf 3.} Consider the representations $V_1, V_2$ and $V_3$ defined before of \eqref{tabaco1}, together with $$V_0 : (t,x,y) \mapsto (\omega_p, 1, 1).$$As before, let $W_j$ be the rational irreducible representation induced by $V_j$. By proceeding as in the previous cases,  the group algebra decomposition of $JC_j$ with respect to $G$ is \begin{equation}\label{raton}JC_j \sim B_{W_0} \times B_{W_1} \times B_{W_2} \times B_{W_3}\end{equation}where each factor has dimension $\tfrac{p-1}{2},$ for each $j.$ The fact that $\rho_{\langle x, y \rangle} \cong W_0$ and $\rho_{\langle x \rangle} \cong W_0 \oplus W_1$ coupled with \cite[Corollary 5.6]{CR06} imply that $$B_{W_0} \sim J(C_j)_{\langle x,y \rangle} \mbox{ and } B_{W_1} \sim \mbox{Prym}((C_j)_{\langle x \rangle} \to (C_j)_{\langle x, y \rangle}).$$Similarly, one obtains that $$B_{W_2} \sim \mbox{Prym}((C_j)_{\langle y \rangle} \to (C_j)_{\langle x, y \rangle}) \mbox{ and } B_{W_3} \sim \mbox{Prym}((C_j)_{\langle xy \rangle} \to (C_j)_{\langle x, y \rangle})$$and the proof follows after replacing the isogenies above in \eqref{raton}.

\s

{\bf 4.} Consider $G=\langle a, b : a^p=b^4=1, [a,b]=1\rangle$ and its complex irreducible representations $$V_1: (a,b) \mapsto (\omega_p, 1), \,\, V_2: (a,b) \mapsto (\omega_p, -1) \, \mbox{ and } \, V_3: (a,b) \mapsto (\omega_p, \sqrt{-1}).$$If we denote the corresponding rational irreducible representations by  $W_1, W_2$ and $W_3$, then the group algebra decomposition of $JS_j$ with respect to $G$ is $$JS_j \sim B_{W_1} \times B_{W_2} \times B_{W_3}$$where the dimension of the first two factors is $\tfrac{p-1}{2}$ whereas the dimension of the latter is $p-1,$ for each $j.$ The proof follows after observing that$$J(S_j)_{\langle b^2 \rangle} \sim B_{W_1} \times B_{W_2}.$$

\subsection*{Proof of Theorem \ref{jaco}}

We start by providing the group algebra decomposition of $JX_p$ with respect to its full automorphism group. The complex irreducible representations of $$G=\langle a, r, s : a^p=r^4=s^2=(sr)^2= [a,r]=[a,s]=1\rangle \cong \mathbb{Z}_p \times \mathbf{D}_4$$ are given by
$${V}_{i j} := \chi_i \otimes \vartheta_j \mbox{ with } 0 \leqslant i \leqslant p-1 \mbox{ and } 0 \leqslant j \leqslant 4,$$
where $\chi_i(a) = \omega_p^i$ and
$$\vartheta_0: (r,s) \mapsto (1,1), \,\, \vartheta_1: (r,s)\mapsto (1,-1), \,\, \vartheta_2: (r,s) \mapsto (-1,1), \,\, \vartheta_3: (r,s)\mapsto (-1,-1)$$ $$\vartheta_4: r \mapsto \begin{psmallmatrix}\sqrt{-1}& 0\\0& -\sqrt{-1}\end{psmallmatrix}, s \mapsto \begin{psmallmatrix}0 & 1\\1 & 0\end{psmallmatrix}$$

\s

In addition, the rational irreducible representations of $ G$ are given by
$${W}_j := {V}_{0 j} = \chi_0 \otimes \vartheta_j  \, \mbox{ and } \,{U}_k := \oplus_{i=1}^{p-1}{V}_{i k}$$where $0 \leqslant j \leqslant 4$ and $0 \leqslant k \leqslant 4.$ It follows that the  group algebra decomposition of $JX_p$ with respect to $G$ is$$JX \sim {B}_{{W}_0} \times   {B}_{{W}_1} \times  {B}_{{W}_2} \times  {B}_{{W}_3} \times  {B}_{{W}_4}^2 \times  {B}_{{U}_0} \times   {B}_{{U}_1} \times  {B}_{{U}_2} \times  {B}_{{U}_3} \times  {B}_{{U}_4}^2.$$A routine computation shows that $\rho_{\langle asr \rangle}= {W}_0 \oplus {W}_3 \oplus {W}_4$, $\rho_{\langle ar \rangle}={W}_0 \oplus{W}_1$ and \begin{equation}\label{rhos}\rho_{\langle s \rangle}={W}_0 \oplus {W}_2 \oplus {W}_4 \oplus {U}_0 \oplus {U}_2 \oplus {U}_4,\end{equation}In this way, by \cite{R07}, we obtain that $\dim {B}_{{W}_j}=0$ for each $0 \leqslant j \leqslant 4$ and \begin{displaymath}
\dim {B}_{{U}_k} = \left\{ \begin{array}{cl}
 \tfrac{p-1}{2}& \textrm{if $k=1,3,4$}\\
 0 & \textrm{if $k=0,2$}
  \end{array} \right.
\end{displaymath}It follows that $JX_p \sim    {B}_{{U}_1} \times  {B}_{{U}_3}  \times  {B}_{{U}_4}^2.$ We now consider the following subgroups of $G$ $$H_1=\langle r \rangle, \,\, H_2=\langle r^2, sr \rangle \,\, \mbox{ and }\,\, H_3=\langle s \rangle,$$and denote by $\mathcal{Z}_j$  the quotient Riemann surface $X_p/H_j$ for $j=1,2,3.$ After considering \cite[Proposition 5.2]{CR06}, the fact that $$\rho_{H_1} =   {W}_0 \oplus {W}_1  \oplus {U}_0  \oplus {U}_1 \mbox{ and } \rho_{H_2}  =   {W}_0 \oplus {W}_3  \oplus {U}_0  \oplus {U}_3$$coupled with \eqref{rhos} imply that  
$$ J\mathcal{Z}_1 \sim   {B}_{{U}_1}, \,\,\,  J\mathcal{Z}_2 \sim  {B}_{{U}_3} \,\, \mbox{ and } \,\, J\mathcal{Z}_3 \sim  {B}_{{U}_4}.$$As a result we obtain that 
\begin{equation}\label{kk}JX_p \sim    J\mathcal{Z}_1 \times J\mathcal{Z}_2  \times  J\mathcal{Z}_3^2,\end{equation}and the genus of $\mathcal{Z}_j$ is $\tfrac{p-1}{2}$ for each $j.$ The fact that  
${\bf N}_G(H_j)/H_j  \cong {\mathbb Z}_{2p}$ for $j=1,2,3$ implies that  $\mathcal{Z}_j$ is endowed with a group of automorphism isomorphic to ${\mathbb Z}_{2p}.$ Thus, $$\mathcal{Z}_1 \cong \mathcal{Z}_2 \cong \mathcal{Z}_3 \cong \hat{X}_p,$$where $\hat{X}_p$ is the unique hyperelliptic among the Riemann surfaces of genus $\tfrac{p-1}{2}$ endowed with an automorphism of prime order $p$ (such Riemann surfaces are usually called {\it Lefschtez surfaces}). It then follows that \begin{equation}\label{mo}JX_p \sim J\hat{X}_p^4.\end{equation}Note that $\hat{X}_p$ is represented by the curve $y^p=x(x-1)^{p-2}$ and the automorphisms of order $2p$  by $(x,y) \mapsto (\tfrac{1}{x},-\omega_p \tfrac{y}{x}).$ Finally, as proved by Lefschetz himself in \cite{L}, $J\hat{X}_p$ is simple, showing that the isogeny \eqref{mo} is the Poincar\'e isogeny decomposition of $JX_p.$

\subsection*{Proof of Theorem \ref{CM}} \mbox{}
We will prove only the first item, as the proof of the remaining ones is analogous. We keep the same notation employed in the proof of Theorem \ref{jaco}.

\s

 We recall that, according to 
\cite[Proposition 2]{br} (see also \cite{e}), the character $\chi$ of the representation of $G$ on 
$H^1(X_p,\mathbb{Q})$ satisfies $$(\rho_{\langle asr \rangle} \oplus \rho_{\langle ar \rangle} \oplus \rho_{\langle s \rangle}) \oplus \chi \cong 2 W_0 \oplus \rho_{\mbox{\tiny regular}}$$and therefore
is given by$$ \chi \cong {U}_1 \oplus  {U}_3 \oplus  {U}_4.
$$Further, according to   \cite[p. 202, Corollaire]{g}, there is
an isomorphism of ${\mathbb Q}[G]$-modules
$$
H^1(X, {\mathbb Q}) \simeq  {U}_1 \oplus  {U}_3 \oplus  {U}_4
$$
and  an isomorphism of ${\mathbb Q}[H_j \backslash G/H_j]$-modules
\begin{equation*} \label{eq2.2}
H^1(\mathcal{Z}_j, {\mathbb Q}) \simeq H^1(X_p, {\mathbb Q})^{H_j} \simeq ({U}_1 \oplus  {U}_3 \oplus  {U}_4)^{H_j} = \left\{ \begin{array}{cr} {U}_1  & \mbox{ for } \; j = 1\\ {U}_3  & \mbox{ for } \; j = 2\\ 
{U}_4  & \mbox{ for } \; j = 3\\  \end{array}\right.
\end{equation*}

This implies that the canonical homomorphism 
${\mathbb Q}[H_j
\backslash G/H_j] \rightarrow \mbox{End}_{{\mathbb Q}}(J\mathcal{Z}_j)$ induces a homomorphism $$\rho_j: {\mathbb Q}[H_j
\backslash G/H_j] \rightarrow \mbox{End}(H^1(\mathcal{Z}_j, {\mathbb Q})) = \mbox{End}( (U_1 \oplus  {U}_3 \oplus  {U}_4)^{H_j}),$$ and the image of
$\rho_j$ is isomorphic to the image of ${\mathbb Q}[H_j \backslash G/H_j]$ in
$\mbox{End}_{{\mathbb Q}}(J\mathcal{Z}_j)$. Now the image of $ {\mathbb Q}[G]$ in $\mbox{End}({U}_1 \oplus  {U}_3 \oplus  {U}_4)$ is isomorphic to  $${\mathbb Q}(\omega_p) \times {\mathbb Q}(\omega_p)
 \times {M}_2(\mathbb{Q}(\omega_p))$$and$${\mathbb Q}[H_1 \backslash G/H_1] \: \cong \: {\mathbb Q} \times {\mathbb Q} \times \mathbb{Q}(\omega_p)  \times  \mathbb{Q}(\omega_p)
$$ $$
{\mathbb Q}[H_2 \backslash G/H_2] \cong {\mathbb Q} \times {\mathbb Q} \times    {\mathbb Q}(\omega_p)  \times  {\mathbb Q}(\omega_p)
$$ $$
{\mathbb Q}[H_3 \backslash G/H_3] \cong  {\mathbb Q} \times {\mathbb Q} \times  {\mathbb Q} \times  {\mathbb Q}(\omega_p)  \times  {\mathbb Q}(\omega_p) \times {\mathbb Q}(\omega_p)
$$ (see \cite[Theorem 4.4]{CR06} and \cite{e}).
Hence, the image of ${\mathbb Q}[H_j \backslash G/H_j]$ in $\mbox{End}( ({U}_1 \oplus  {U}_3 \oplus  {U}_4)^{H_j}) $ is isomorphic to ${\mathbb Q}(\omega_p),$
and therefore
$${\mathbb Q}(\omega_p) \subset \mbox{End}_{{\mathbb Q}}(J\mathcal{Z}_j) \mbox{ for }j=1,2,3.$$The result follows from the isogeny \eqref{kk} given in the proof of Theorem \ref{jaco}.

\end{document}